\documentclass[onefignum,onetabnum]{siamart250211}

\usepackage{theorem}
\usepackage{enumerate}
\usepackage{array}
\usepackage{amssymb}
\usepackage{mathtools}
\usepackage{latexsym}
\usepackage{makeidx}
\usepackage{fancybox}
\input epsf.sty
\usepackage{subcaption}

\newtheorem{example}[theorem]{Example}

\numberwithin{equation}{section}

\outer\def\proclaim #1. #2\par{\medbreak \noindent{\bf#1.\enspace}{\sl#2}\par
  \ifdim\lastskip<\medskipamount
  \removelastskip\penalty55\medskip\fi}
\def\state #1. { \noindent{\bf#1.\enspace}}
\def\algo #1. { \noindent{\bf#1.\enspace}}

\DeclareMathOperator*{\argmin}{argmin}

\DeclareMathOperator{\cl}{cl}
\DeclareMathOperator{\bdry}{bdry}

\DeclareMathOperator{\exs}{exs}
\DeclareMathOperator{\dom}{dom}

\DeclareMathOperator{\epi}{epi}
\DeclareMathOperator{\hypo}{hypo}
\DeclareMathOperator{\nt}{int}

\newcommand{\comp}{\,{\raise 1pt \hbox{$\scriptstyle\circ$}}\,}

\newcommand{\reals}{\mathbb{R}}
\newcommand{\Reals}{\overline{\mathbb{R}}}
\newcommand{\natnums}{{{\rm l} \kern -.13em {\rm N} }}
\newcommand{\nats}{\mathbb{N}}
\newcommand{\snats}{{I\kern -.29em N}}
\newcommand{\rats}{{Q\kern -.64em \raise 1pt \hbox{$\scriptstyle |$}\;\,}}
\newcommand{\srats}
	{{Q\kern -.56em \raise 1.2pt \hbox{$\scriptscriptstyle /$}\,}}
\newcommand{\ints}{Z\kern -.46em Z}
\newcommand{\ball}{\mathbb{B}}
\newcommand{\pluss}{\hskip1pt \raise1pt\vbox{\hrule width6pt \vskip1pt \hrule
                    width6pt} \kern-4pt{\lower1pt\hbox{\vrule height6pt
		    \kern1pt\vrule height6pt}}\hskip5pt}
\newcommand{\eop}
	{\hfill{$\vcenter{\hrule height1pt \hbox{\vrule width1pt height5pt
   	 \kern5pt \vrule width1pt} \hrule height1pt}$} \medskip}
\newcommand{\half}
	{{\raisebox{1pt}{$\frac{1}{2}$}}}

\newcommand{\setd}{{ d \kern -.15em l}}
\newcommand{\hatsetd}{ d \hat{\kern -.15em l }}

\renewcommand{\epsilon}{\varepsilon}
\renewcommand{\phi}{\varphi}

\newcommand{\tto}{\;{\lower 1pt \hbox{$\rightarrow$}}\kern -12pt
           \hbox{\raise 2.5pt \hbox{$\rightarrow$}}\;}
\newcommand{\overto}[1]{\,{\raise 0pt\hbox{$\rightarrow$}}\kern -9pt
     \hbox{\lower 3pt \hbox{$\scriptscriptstyle#1$}}\hskip6pt}
\newcommand{\underto}[1]{\,{\lower 1pt\hbox{$\rightarrow$}}\kern -9pt
     \hbox{\raise 4pt \hbox{$\,\scriptscriptstyle#1$}}\hskip7pt}
\newcommand{\bigoverto}[1]{{\raise 0pt\hbox{$\,\longrightarrow$}}\kern -16pt
     \hbox{\lower 3pt \hbox{$\scriptscriptstyle#1$}}\hskip4pt}
\newcommand{\bigunderto}[1]{\,{\lower 1pt\hbox{$\longrightarrow$}}\kern -16pt
     \hbox{\raise 4pt \hbox{$\,\scriptscriptstyle#1$}}\hskip6pt}
\newcommand{\bigbigto}[2]{\,{\raise 0pt\hbox{$\,\longrightarrow$}}\kern -16pt
     \hbox{\lower 3pt \hbox{$\scriptscriptstyle#2$}}\kern -10pt
     \hbox{\raise 4pt \hbox{$\,\scriptscriptstyle#1$}}\hskip7pt}
\newcommand{\downto}{{\raise 1pt \hbox{$\scriptscriptstyle \,\searrow\,$}}}
\newcommand{\upto}{{\raise 1pt \hbox{$\scriptscriptstyle \,\nearrow\,$}}}

\newcommand{\notimply}
	{\quad\hbox{$\Longrightarrow \kern -14pt {/}$}\hskip6pt\quad}

\newcommand{\lto}{\,{\lower 1pt\hbox{$\rightarrow$}}\kern -10pt
     \hbox{\raise 4pt \hbox{$\, \scriptstyle l$}}\hskip7pt}
\newcommand{\eto}{\,{\lower 1pt\hbox{$\rightarrow$}}\kern -10pt
     \hbox{\raise 4pt \hbox{$\, \scriptstyle e$}}\hskip7pt}
\newcommand{\hto}{\,{\lower 1pt\hbox{$\rightarrow$}}\kern -11pt
     \hbox{\raise 4pt \hbox{$\, \scriptstyle h$}}\hskip7pt}
\newcommand{\pto}{\,{\lower 1pt\hbox{$\rightarrow$}}\kern -11pt
     \hbox{\raise 4.5pt \hbox{$\, \scriptstyle p$}}\hskip7pt}
\newcommand{\cto}{\,{\lower 1pt\hbox{$\rightarrow$}}\kern -11pt
     \hbox{\raise 4pt \hbox{$\, \scriptstyle c$}}\hskip7pt}
\newcommand{\gto}{\,{\lower 1pt\hbox{$\rightarrow$}}\kern -11pt
     \hbox{\raise 4.5pt \hbox{$\, \scriptstyle g$}}\hskip7pt}
\newcommand{\sto}{\,{\lower 1pt\hbox{$\rightarrow$}}\kern -11pt
     \hbox{\raise 4pt \hbox{$\, \scriptstyle s$}}\hskip7pt}
\newcommand{\awto}{\,{\lower 1pt\hbox{$\rightarrow$}}\kern -15pt
     \hbox{\raise 4pt \hbox{$\, \scriptstyle aw$}}\hskip7pt}
\def\Nto{\,{\raise 1pt\hbox{$\rightarrow$}}\kern -12pt
     \hbox{\lower 3pt \hbox{$\, \scriptstyle N$}}\hskip7pt}
\def\Cto{\,{\raise 1pt\hbox{$\rightarrow$}}\kern -14pt
     \hbox{\lower 3pt \hbox{$\, \scriptstyle C$}}\hskip7pt}
\def\fto{\,{\raise 1pt\hbox{$\rightarrow$}}\kern -14pt
     \hbox{\lower 3pt \hbox{$\, \scriptstyle f$}}\hskip7pt}

\newcommand{\low}[1]{{\lower1pt \hbox{$\scriptstyle #1$}}}
\newcommand{\loww}[1]{{\lower2pt \hbox{$\scriptstyle #1$}}}
\newcommand{\high}[1]{{\raise1pt \hbox{$\scriptstyle #1$}}}

\newcommand{\nsum}{\mathop{\sum}\nolimits}

\newcommand{\nliminf}{\mathop{\rm liminf}\nolimits}

\newcommand{\nlimsup}{\mathop{\rm limsup}\nolimits}
\newcommand{\ninf}{\mathop{\rm inf}\nolimits}
\newcommand{\nsup}{\mathop{\rm sup}\nolimits}

\newcommand{\nnmin}{\mathop{\rm minimize}}
\newcommand{\nnmax}{\mathop{\rm maximize}}

\newcommand{\nargmin}{\mathop{\rm argmin}\nolimits}

\newcommand{\lwdy}[2]{\mathrel{\mathop
        {\raisebox{0.1ex}{\null$#1$}}{\hbox{\kern -1.0em
	{\raisebox{-0.8ex}{$\scriptstyle{\;\to #2}$}}}}}}
\newcommand{\lwwdy}[2]{\mathrel{\mathop
        {\raisebox{0.2ex}{\null$#1$}}{\hbox{\kern -1.0em
	{\raisebox{-1.1ex}{$\scriptstyle{\;\to #2}$}}}}}}
\newcommand{\slwdy}[2]{\scriptsize{{\mathrel{\mathop
        {\raisebox{0.1ex}{\null$#1$}}{\hbox{\kern -1.0em
	{\raisebox{-0.8ex}{$\scriptstyle{\;\to #2}$}}}}}}}}
\newcommand{\slwwdy}[2]{\scriptsize{{\mathrel{\mathop
        {\raisebox{0.2ex}{\null$#1$}}{\hbox{\kern -1.0em
	{\raisebox{-1.1ex}{$\scriptstyle{\;\to #2}$}}}}}}}}

\definecolor{lightgray}{gray}{0.75}
\definecolor{myred}{rgb}{0.55,0,0}
\definecolor{myblue}{rgb}{0,0,0.5} 
\definecolor{mygreen}{rgb}{0,0.5,0} 
\definecolor{purple}{rgb}{0.5,0,0.5} 
\definecolor{turq}{rgb}{0,0.805,0.816} 
\definecolor{maroon}{rgb}{0.51,0,0}
\definecolor{MAROON}{rgb}{0.51,0,0}
\definecolor{redor}{rgb}{0.78,0.078,0.078}
\definecolor{dgreen}{rgb}{0,0.3,0}

\newcommand{\bcdot}{\,{\raise .2ex \hbox{$\centerdot$}}\,}

\newcommand{\beq}{\begin{equation}}
\newcommand{\eeq}{\end{equation}}
\newcommand{\ba}{\begin{array}}
\newcommand{\ea}{\end{array}}

\newcommand{\redrev}[1]{{\color{black}{#1}}}

\newcommand{\redrevvv}[1]{{\color{black}{#1}}}

\ifpdf
\hypersetup{
	pdftitle={Approximations of Rockafellians, Lagrangians, and Dual Functions},
	pdfauthor={J. Deride, and J.O. Royset}
}
\fi

\title{Approximations of Rockafellians, Lagrangians, and Dual Functions}

\author{
	Julio Deride
	\thanks{Faculty of Engineering and Science, Universidad Adolfo Ib\'a\~nez, Santiago, Chile. Email: \email{julio.deride@uai.cl}.}
\and Johannes O. Royset
\thanks{Daniel J. Epstein Department of Industrial \& Systems Engineering, University of Southern California, Los Angeles, CA 90089. Email: \email{royset@usc.edu}.}
}

\begin{document}
	
	\maketitle
	
	\begin{abstract}
Solutions of an optimization problem are sensitive to changes caused by approximations or parametric perturbations, especially in the nonconvex setting. This paper shows that solutions of substitute problems, constructed from Rockafellian functions, can be less sensitive to such changes. Unlike classical stability analysis focused on local changes around (local) minimizers, we employ epi-convergence to examine whether approximating or perturbed problems suitably approach an actual (unperturbed) problem globally. \redrevvv{We demonstrate that solutions derived from the Rockafellian-based substitute problems converge to solutions of the actual optimization problem under suitable conditions, providing a rigorous alternative to potentially unstable direct approximations.} We quantify the rates of convergence that often lead to Lipschitz-kind stability properties for the substitute problems.
	\end{abstract}
	


\section{Introduction}\label{sec:intro}

Given the problem of minimizing $\phi:\reals^n\to \Reals = [-\infty,\infty]$, we are often compelled to consider one or more approximating problems. The approximations may stem from imprecise implementation of conceptual operations such as integration of nontrivial functions and maximization over uncountable sets or from smoothing methods and penalty functions. Approximations also arise in sensitivity analysis as model parameters are perturbed, for example in an effort to identify the problem's robustness to data manipulation. These approximations may cause disproportionally large changes in minimizers, minimum values, and also level-sets; see, e.g., \cite{Royset.20b,RoysetChenEckstrand.22}. Fortunately, every optimization problem can be associated with many {\em substitute problems} and they may possess more desirable stability properties. This paper develops sufficient conditions for epigraphical and hypographical stability of substitute problems defined by Rockafellian relaxations, Lagrangian relaxations, and dual functions. 

Stability of optimization problems can be viewed from different angles; see, for example, the monographs \cite{Attouch.84,AubinEkeland.84,VaAn,BonnansShapiro.00,Mordukhovich.18}. Metric regularity and calmness give rise to local stability \cite{IoffeOutrata.08,Penot.10}. Tilt-stability
\cite{PoliquinRockafellar.98,EberhardWenczel.12,LewisZhang.13,DrusvyatskiyLewis.13} and full-stability
\cite{MordukhovichRockafellarSarabi.13,MordukhovichNghiaRockafellar.15} endow local minimizers with Lipschitz-type properties under perturbations and are supported by an extensive calculus for coderivatives; see \cite{BenkoRockafellar.24} for recent efforts in this direction. An alternative approach to Lipschitz stability based on a growth condition appears in \cite{Klatte.96}. Directional regularity may also provide quantitative estimates of the rate of change of a min-value function \cite{Gollan.84,BonnansCominetti.96,BonnansShapiro.00}. For specific developments in stochastic programming, we mention \cite{HenrionRoemisch.04,DentchevaHenrionRuszczynski.07}.

Aligned with the pioneering work in \cite{AttouchWets.91,AttouchWets.93a,AttouchWets.93b}, this paper takes an alternative, global perspective afforded by epi-convergence and its quantification in terms of the truncated Hausdorff distance. With $\phi^\nu:\reals^n\to \Reals$, $\nu\in \nats=\{1, 2, \dots \}$ being a sequence of approximating functions, a concerning situation arises when $\phi^\nu$ fails to epi-converge to the actual function $\phi$ as $\nu\to \infty$; we include several examples in the following sections. In the absence of this basic convergence property, one cannot expect minimizers, minimum values, and level-sets of the approximations to approach those of the actual function; see \cite[Chap. 7]{VaAn}. Nevertheless, Rockafellian relaxations, Lagrangian relaxations, and dual problems can be brought in and they often achieve the desirable convergence of near-minimizers, minimum values, and level-sets at quantifiable rates. \redrevvv{Under broad conditions involving epi-convergence and exactness, cluster points of minimizers from these relaxations are, in fact, minimizers of the actual underlying problem, a property not generally shared by direct approximations.}

While stability results near a local solution leverage concepts from second-order subdifferentiability, supporting calculus rules, and notions of ``regularity'' (cf. \cite{VaAn,BonnansShapiro.00,Mordukhovich.18}), our perspective furnishes different kinds of insights using different kinds of tools. We measure the discrepancy between two optimization problems using the epigraphs or hypographs of the corresponding functions, with subdifferentiability playing a secondary role. We even allow for variables with integer restrictions. The approach is advantageous in the sense that each level-set of one of the functions can be related to a level-set of the other function, which amounts to a global quantification of the discrepancy between the problems; see \cite[Prop. 2.3]{Royset.20b} and below for details. We also avoid the typical regularity conditions such as constraint qualifications, which are central to any local stability analysis \cite{VaAn,BonnansShapiro.00,Mordukhovich.18}. Instead tightness enters as seen below. Despite some calculus rules \cite{AzePenot.90b,AttouchWets.91,Royset.20b}, our approach comes with the challenge of computing truncated Hausdorff distances between epigraphs and also between hypographs. A contribution of this paper is to demonstrate that such computations are indeed possible in practically important settings. 

All the substitute problems under consideration stem from Rockafellian functions. By definition, $f:\reals^m\times\reals^n\to \Reals$ is a {\em Rockafellian} for $g:\reals^n\to\Reals$ when $f(0,x) = g(x)$ for all $x\in\reals^n$ \cite[Sec. 5.A]{primer}. As seen below, epi-convergence of Rockafellians is a key property in the analysis of the various substitute problems and their approximations. This perspective can be traced back to \cite{Bergstrom.80}, but was extended in \cite{AttouchAzeWets.86,Attouch.86,AtAW88:sad,Aze.88,AzeRahmouni.96}, especially to infinite-dimensional spaces, using epi/hypo-convergence \cite{AttouchWets.83}. These references deal with convex Rockafellians and the resulting convex-concave Lagrangians. In the present paper, we examine Rockafellians in the {\em nonconvex} setting and also quantify the rate of convergence. 

Suppose that $\{f,f^\nu:\reals^m\times\reals^n\to \Reals, \nu\in\nats\}$ are Rockafellians for $\{\phi, \phi^\nu, \nu\in \nats\}$, respectively. They define the following substitute problems as alternatives to the {\em actual problem} of minimizing $\phi$ and the {\em approximating problem} of minimizing $\phi^\nu$; the notation applies throughout the paper:

The {\em Rockafellian relaxations} (given $y,y^\nu\in\reals^m$)
\begin{align*}
  &\nnmin_{u\in\reals^m, x\in \reals^n} f_y(u,x), && \mbox{ where } f_y(u,x) = f(u,x) - \langle y, u\rangle\\
  &\nnmin_{u\in\reals^m, x\in \reals^n} f^\nu_{y^\nu}(u,x), && \mbox{ where } f^\nu_{y}(u,x) = f^\nu(u,x) - \langle y, u\rangle.
\end{align*}

The {\em Lagrangian relaxations} (given $y,y^\nu\in\reals^m$)
\begin{align*}
  &\nnmin_{x\in \reals^n} \, l(x,y), &&\mbox{ where } l(x,y) = \ninf_u \big\{f(u,x) - \langle y, u\rangle\big\}\\
  &\nnmin_{x\in \reals^n} \, l^\nu(x,y^\nu), &&\mbox{ where } l^\nu(x,y) = \ninf_u \big\{f^\nu(u,x) - \langle y, u\rangle\big\}.
\end{align*}

The {\em dual problems}
\begin{align*}
  &\nnmax_{y\in \reals^m} \, \psi(y), &&\mbox{ where } \psi(y) = \ninf_{u,x} \big\{f(u,x) - \langle y, u\rangle\big\}\\
  &\nnmax_{y\in \reals^m} \, \psi^\nu(y), &&\mbox{ where } \psi^\nu(y) = \ninf_{u,x} \big\{f^\nu(u,x) - \langle y, u\rangle\big\}.
\end{align*}

Since a Rockafellian simply parametrizes the underlying function, there are countless possible versions of these substitute problems with concrete examples following below. The concept of developing dual problems through parametrization was pioneered in \cite[Chap. 29]{Rockafellar.70}, with further developments in \cite{Rockafellar.74,Rockafellar.85}, but the idea of embedding a problem within a parametric family extends back to \cite{Rockafellar.63}. The name ``Rockafellian'' emerges in \cite[Chap. 5]{Royset.21,primer}, with ``bifunction,'' ``bivariate function,'' and  ``perturbation function'' also being used in the literature. For common choices of Rockafellians, the importance of Lagrangian relaxations and dual problems is well documented (see, e.g., \cite[Chap. 5]{primer}).

We refer to $f_y$ and $f^\nu_{y^\nu}$ as {\em tilted Rockafellians} due to their construction from $f$ and $f^\nu$ via a linear term. The bivariate functions $l$ and $l^\nu$ are called {\em Lagrangians}, with classical forms from nonlinear programming arising for special choices of Rockafellians; see Example \ref{e:composite} below.  We refer to $y,y^\nu$ in Rockafellian and Lagrangian relaxations as {\em multiplier vectors}. The optimization of {\em dual functions} $\psi,\psi^\nu$ can be viewed as a process of finding the best multiplier vectors.

\begin{figure}
\centering
\includegraphics[width=0.83\textwidth]{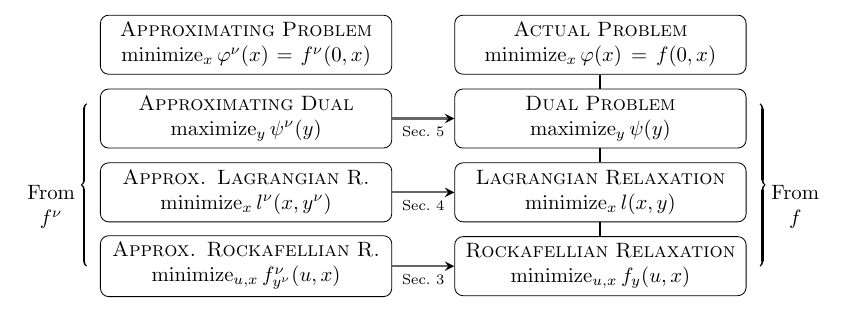}
\caption{Overview of problems defined by Rockafellians $f,f^\nu$.}\label{overview}
\end{figure}

Figure \ref{overview} provides a road map for the paper. The absence of an arrow between the problems on top indicates that the approximating functions $\phi^\nu$ may fail to epi-converge to $\phi$. Nevertheless, a Rockafellian $f$ for $\phi$ provides substitute problems listed on the right in Figure \ref{overview} that might be better behaved in the sense that if approximations are introduced in $f$ to produce $f^\nu$, a Rockafellian for $\psi^\nu$, then the resulting tilted Rockafellians $f_{y^\nu}^\nu$ epi-converge to $f_y$ and the Lagrangians $l^\nu(\,\cdot\,, y^\nu)$ epi-converge to $l(\,\cdot\,, y)$ when $y^\nu\to y$. This is represented by the two lowest arrows in Figure \ref{overview}. Likewise, we aspire to achieve hypo-convergence of the dual functions $\psi^\nu$ to $\psi$ as illustrated by an arrow in the figure. (The shift to hypo-convergence is dictated by the convention of stating dual problems in terms of maximization instead of minimization.) Without relying on convexity, we show that the convergences indicated by the horizontal arrows in the figure follow naturally when $f^\nu$ epi-converges to $f$ and that the rates can be quantified in a nonasymptotic analysis. These results appear to be the first of their kind. The prior effort in \cite{RoysetChenEckstrand.22} limits the scope to certain stochastic programs and omit a discussion of Lagrangian relaxation and dual problems. The papers \cite{BurkeHoheisel.13,Royset.22a} give sufficient conditions for epi-convergence of approximating functions in the context of composite problems, but neither address substitute problems nor quantification of the convergence. In fact, the delicate assumptions adopted in \cite[Thm. 2.4]{Royset.22a} to achieve epi-convergence highlight the need for passing to substitute problems. While discussing epi-convergence and its quantification broadly, \cite{Royset.18,Royset.20b} fail to address the specific problem structures appearing from Rockafellian relaxation, Lagrangian relaxation, and dual problems. We also provide a novel direction to study dual problems that is distinct from those based on lopsided convergence \cite{AttouchWets.83b,RoysetWets.19a}. 

Even if the substitute problems on the left converge to those on the right in Figure \ref{overview}, the fact remains that the substitute problems on the right are mere relaxations of the actual problem. However, we can bridge that gap using an extensive theory about strong duality (e.g., \cite[Thms. 11.39, 11.59]{VaAn}), adjusted for our purpose in Definition \ref{def:exactness} using the name {\em exactness}, which follows \cite{RoysetChenEckstrand.22}; see also \cite[Def. 2.146]{BonnansShapiro.00}. The definition posit conditions for ``equivalence'' among all four problems on the right-hand side in the figure as indicated by short vertical lines. Thus, we envision the analysis for an application to proceed in two steps: First, identify a Rockafellian $f$ such that equivalence is ensured on the right. Second, introduce the approximations to obtain $f^\nu$ and then confirm one or more of the ``convergence arrows" from left to right. Even if the first step is not achieved, the substitute problems on the right-hand side are relaxations of the actual problem, and convergence to such relaxations by the approximations on the left-hand side can still be meaningful. 

The paper starts with preliminaries in Section 2. Sections 3, 4, and 5 discuss tilted Rockafellians, Lagrangians, and dual functions, respectively. Section 6 provides supplementary results.

\section{Preliminaries} 

We follow the notation of \cite{primer}. For $C\subset \reals^n$, $\iota_C(x) = 0$ if $x\in C$ and $\iota_C(x) = \infty$ otherwise. We write $\nt C$ for the {\em interior} of $C$, $\cl C$ for the {\em closure} of $C$, and $\bdry C=\cl C\setminus \nt C$ for the {\em boundary} of $C$. We write $C^\nu\sto C$ when $C^\nu\subset \reals^n$ {\em set-converges}  to $C$ in the sense of Painlev\'{e}-Kuratowski. If $C$ is a convex set, then its {\em normal cone} at $x$ is denoted by $N_C(x)$.

For a subsequence $N\subset \nats$, we write $x^\nu \Nto x$ when $\{x^\nu, \nu\in N\}$ converges to $x$. We adopt the convention $\infty - \infty = \infty$; see also \cite[Sec. 1.D]{primer}. Inequalities between vectors are interpreted componentwise.

For a function $g:\reals^n\to \Reals$, its {\em domain} is $\dom g = \{x\in \reals^n\,|\,g(x)<\infty\}$ and $\alpha$-{\em level-set} is $\{g \leq \alpha\} = \{x\in\reals^n \, | \, g(x) \leq \alpha\}$. Its {\em epigraph} is $\epi g = \{(x,\alpha) \in \reals^n\times \reals \,| \, g(x) \leq \alpha\}$. The {\em minimum value} of $g$ is $\inf g = \inf_x g(x) = \inf\{g(x)\,|\,x\in\reals^n\}$ and the set of {\em minimizers} is $\nargmin g = \nargmin_x g(x) = \{x^\star \in \dom g\,|\,g(x^\star) \leq \inf g\}$. For any $\epsilon \in (0,\infty)$, the set of {\em near-minimizers} is written as $\epsilon\mbox{-}\nargmin g = \{x^\star \in \dom g\,|\,g(x^\star) \leq \inf g + \epsilon\}$. We say that $g$ is {\em proper} when $g(x)>-\infty$ for all $x\in \reals^n$ and $\epi g \neq \emptyset$. It is {\em lower semicontinuous} (lsc) when $\epi g$ is closed and it is {\em convex} if $\epi g$ is convex. We denote by $g^*$ the {\em conjugate} of $g$, i.e., $g^*(y) = \sup_x \langle x,y\rangle - g(x)$. If $g$ is convex, its set of {\em subgradients} at $\bar x\in \reals^n$ is written as $\partial g(\bar x)$. 
We say $g$ is {\em Lipschitz with modulus function} $\kappa$ when it is real-valued, $\kappa:[0,\infty)\to[0,\infty)$, and $|g(x) - g(\bar x)| \leq \kappa(\rho)\|x - \bar x\|_2$ whenever $\|x\|_2\leq \rho$, $\|\bar x\|_2\leq \rho$, and  $\rho\in [0,\infty)$. 

The functions $g^\nu:\reals^n\to \Reals$ {\em epi-converge} to $g:\reals^n\to \Reals$, written $g^\nu\eto g$, when
\begin{align}
  &\forall x^\nu \to x, ~~\nliminf g^\nu(x^\nu) \geq g(x),\label{eqn:liminf}\\
  &\forall x, ~\exists x^\nu\to x \mbox{ with } \nlimsup g^\nu(x^\nu)\leq g(x).\label{eqn:limsup}
\end{align}
We say that $\{g^\nu, \nu\in\nats\}$ is {\em tight} if for all $\epsilon>0$, there exist a compact set $B_\epsilon \subset \reals^n$ and an integer $\nu_\epsilon \in \nats$ such that $\ninf \{g^\nu(x)~|~x\in B_\epsilon\} \leq \inf g^\nu + \epsilon$ for all $\nu\geq \nu_\epsilon$. The consequences of epi-convergence for minimizers and minimum values are well known; see Proposition \ref{prop:tigthepi} for a summary.

The functions $g^\nu:\reals^n\to \Reals$ {\em hypo-converge} to $g:\reals^n\to \Reals$, written $g^\nu\hto g$, when $-g^\nu \eto -g$.

The dual function produced by a Rockafellian $f$ is $\psi = -f^*(\,\cdot\,,0)$ by definition, and the resulting Lagrangian $l$, for any $x$, has $-l(x, \,\cdot\,)$ as the conjugate of $f(\,\cdot\,,x)$. An example from the far-reaching area of composite optimization (see, e.g., \cite{CuiPangSen.18,BurkeHoheiselNguyen.21,Royset.22a}) allows for more specific formulas.

\begin{example}{\rm (composite optimization).}\label{e:composite}
For $X,X^\nu\subset\reals^n$, $g_0,g_0^\nu:\reals^n\to \reals$, $G,G^\nu:\reals^n\to \reals^m$, and $h,h^\nu:\reals^m\to \Reals$, suppose that the actual and approximating problems are defined by
\[
\phi(x) = \iota_X(x) + g_0(x) + h\big(G(x)\big) ~~~\mbox{ and } ~~~\phi^\nu(x) = \iota_{X^\nu}(x) + g_0^\nu(x) + h^\nu\big(G^\nu(x)\big).
\]
Even if each ``component'' $X$, $g_0$, $h$, and $G$ of the actual problem is well approximated by $X^\nu$, $g_0^\nu$, $h^\nu$, and $G^\nu$ respectively, $\phi^\nu$ may not epi-converge to $\phi$. For instance, this failure takes place in the simple nonlinear programming setting with $\phi(x) = -x + \iota_{(-\infty,0]}(x^3 - x^2 - x + 1)$ and $\phi^\nu(x) = -x + \iota_{(-\infty,0]}(x^3 - x^2 - x + 1 + 1/\nu)$; see \cite{Royset.20b,Royset.22a} for other problematic instances. 

A possible Rockafellian for $\phi$ and the resulting Lagrangian are given by 
\[
f(u,x) = \iota_X(x) + g_0(x) + h\big(G(x)+u\big)~\mbox{ and }~l(x,y) = \iota_X(x) + g_0(x) + \big\langle G(x), y\big\rangle - h^*(y),
\]
with definitions of the approximating counterparts $f^\nu$ and $l^\nu$ following similarly. While $\phi^\nu$ may not epi-converge to $\phi$, we confirm that the resulting (tilted) Rockafellians epi-converge (Proposition \ref{prop:composite}), the Lagrangians $l^\nu(\,\cdot\,,y^\nu)$ epi-converge to $l(\,\cdot\,, y)$ (Subsection \ref{subsec:compositeLagr}), and the dual functions $\psi^\nu$ hypo-converge to $\psi$ (Subsection \ref{subsec:compositeDual}) under natural assumptions. We also quantify the rate of convergence. 

An explicit formula for the resulting dual function is available, for instance, when $X = \reals^n$, $G(x) = Ax - b$, and $h(y) = \iota_{\{0\}^m}(y)$. Then, $\psi(y) = -\langle b,y\rangle - g_0^*(-A^\top y)$. 
\end{example}

Given a norm $\|\cdot\|$ on $\reals^n$, the {\em centered ball} with radius $\rho \in [0,\infty)$ is $\ball(\rho) = \{ x\in \reals^n | \|x\|\leq \rho \}$. The {\em excess} of $C\subset\reals^n$ over $D\subset\reals^n$, if both sets are nonempty, is $\exs(C; D) = \nsup_{x\in C} \inf_{y\in D}\|x-y\|$, while $\exs(C; \emptyset) = \infty$ for $C\neq \emptyset$ and $\exs(\emptyset; D) = 0$ regardless of $D$. For $\rho\geq 0$, the {\em truncated Hausdorff distance} between $C$ and $D$ is 
\[
\hatsetd_{\rho}(C,D) = \max\big\{\exs( C \cap \ball(\rho); D ), ~\exs( D \cap \ball(\rho); C)\big\}.
\]
The following bounds stem from \cite[Thm. 6.56]{primer} and \cite[Prop. 2.3]{Royset.20b}. They are sharp as discussed in \cite[Sec. 6.J]{primer}. 

\begin{proposition}\label{tapproxoptimalvalue}{\rm (approximation error \cite{primer,Royset.20b}).}
For $g,h:\reals^n\to\Reals$ and $\rho\in [0,\infty)$, adopt any norm $\|\cdot\|$ on $\reals^n$ and the norm $(x,\alpha) \mapsto \max\{\|x\|,|\alpha|\}$ on $\reals^{n+1}$. The following hold: 
\begin{enumerate}[(a)]
    \item If $\epsilon \in [0,2\rho]$, $\delta > \epsilon + 2\hatsetd_{\rho}(\epi g, \epi h)$, the sets $\nargmin g \cap \ball(\rho)$ and $\nargmin h \cap \ball(\rho)$ are nonempty, and $\ninf g, \ninf h \in [-\rho,\rho-\epsilon]$, then 
\begin{align*}
&|\ninf g - \ninf h | \leq \hatsetd_{\rho}(\epi g, \epi h)\\
&\exs\big(\epsilon\mbox{-}\nargmin h \cap \ball(\rho); ~\delta\mbox{-}\nargmin g\big)  \leq \hatsetd_{\rho}(\epi g, \epi h).
\end{align*}

\item If $\epsilon \in [-\rho,\rho]$ and $\delta > \epsilon + \hatsetd_{\rho}(\epi g, \epi h)$, then 
\[
\exs\big(\{h \leq \epsilon\} \cap \ball(\rho); ~\{g \leq \delta\}\big) \leq \hatsetd_{\rho}(\epi g, \epi h).
\]

\end{enumerate}
\end{proposition}

The usefulness of a dual problem depends in part on whether $\sup \psi = \inf \phi$. By Sion's theorem \cite[Cor. 3.3]{Sion.58}, this holds in Example \ref{e:composite} when $X$ is compact, $g_0$ is convex and lsc, $h$ is proper, lsc, and convex, and $\langle G(\cdot),y\rangle$ is convex for all $y\in \dom h^*$; see \cite[Thms. 11.39, 11.59]{VaAn} for other conditions. 

We follow \cite{RoysetChenEckstrand.22} and rely on exactness to confirm the equivalence between the problems on the right-hand side of Figure \ref{overview}. The terminology is motivated by exact penalization \cite{Burke.91}, augmented Lagrangians \cite[Sec. 11.K]{VaAn}, and recent work on linear penalties \cite{Dol16unifying}. It is closely related to \cite[Def. 2.146]{BonnansShapiro.00}, which uses the name ``calmness.'' While \cite[Def. 11.60]{VaAn} identifies exactness with a Lagrangian $l(\,\cdot\,,y)$ having identical minimizers and minimum value to those of $\phi$, we shift the focus to $\inf_x f(\,\cdot\,, x)$. 

\begin{definition}{\rm (exact Rockafellian \cite{RoysetChenEckstrand.22}).}\label{def:exactness}
A Rockafellian $f:\reals^m\times \reals^n\to \Reals$ for $\phi:\reals^n\to\Reals$ is {\em exact, supported by} $y \in \reals^m$, if
\begin{equation}\label{eqn:defexactness}
\ninf_x f(u,x) \geq \inf \phi + \langle y, u\rangle ~~~\forall u\in\reals^m.
\end{equation}
If the inequality holds strictly for all $u\neq 0$, then $f$ is {\em strictly exact}, supported by $y$.
\end{definition}

In terms of the corresponding dual function $\psi$, the condition \eqref{eqn:defexactness} holds if and only if $\psi(y) = \inf \phi$ because generally $\psi(y') \leq \inf \phi$ for all $y'\in \reals^m$ and \eqref{eqn:defexactness} amounts to having $\psi(y) = \inf f_y \geq \inf \phi$. Consequently, when $f$ is exact for $\phi$, supported by $y$, one has 
\begin{equation}\label{eqn:defexactness2}
\inf f_y = \ninf_x l(x,y) = \psi(y) = \sup \psi  = \inf \phi.
\end{equation}
The problems on the right-hand side of Figure \ref{overview} are therefore equivalent in the sense of optimal objective function values; see \cite[Prop. 3.3]{RoysetChenEckstrand.22} and \cite[Sec. 2.5]{BonnansShapiro.00} for further characterizations. Combining this equivalence with epi-convergence, we can formally link Rockafellian and Lagrangian relaxations to the actual problem in the manner indicated in Figure \ref{overview}. For example, we obtain the following fact by using \cite[Prop. 3.3]{RoysetChenEckstrand.22}, Proposition \ref{prop:tigthepi}(d), and Theorem \ref{thm:errorRock}: 
\begin{proposition}{\rm (convergence to actual solution).}\label{prop:convtoactual}
For Rockafellians $\{f,f^\nu:\reals^m \times \reals^n\to \Reals, \nu\in\nats\}$, with $f$ being strictly exact for $\phi$, supported by $y \in \reals^m$, suppose that $f^\nu\eto f$, $y^\nu\to y$, $\epsilon^\nu\to 0$, $\ninf f < \infty$, and  $(\bar u, \bar x)$ is a cluster point of 
\[
\big\{(u^\nu,x^\nu)\in \epsilon^\nu\mbox{-}\nargmin_{u,x} f^\nu_{y^\nu}(u,x), ~\nu\in\nats \big\}. 
\]
Then, $\bar x\in \nargmin \phi$.  
\end{proposition}

Any Rockafellian $f$ can be augmented to produce a strictly exact Rockafellian $(u,x)\mapsto f(u,x) + \iota_{\{0\}^m}(u)$ with other possibilities being discussed below.

\section{Convergence of Tilted Rockafellians}\label{sec:tilted}

The convergence of Rockafellians and the rate carry over to their tilted counterparts nearly without assumptions. This formalizes the bottom arrow in Figure \ref{overview}. A short proof appears in Section \ref{sec:supp}. 

\begin{theorem}{\rm (tilted Rockafellians).}\label{thm:errorRock}
For Rockafellians $\{f,f^\nu:\reals^m \times \reals^n\to \Reals, \nu\in\nats\}$, one has: 

\begin{enumerate}[(a)]

\item If $f^\nu \eto f$ and $y^\nu\to y\in \reals^m$, then $f^\nu_{y^\nu} \eto f_y$.

\item If $\epi f$ and $\epi f^\nu$ are nonempty, $\rho \in [0,\infty)$, $y,y^\nu\in\reals^m$, and $\bar\rho \geq \rho(1 + \max\{\|y\|_2, \|y^\nu\|_2\})$, then
\[
\hatsetd_\rho(\epi f_y, \epi f^\nu_{y^\nu}) \leq \Big( 1 + \max\big\{\|y\|_2, \|y^\nu\|_2\big\} \Big)\hatsetd_{\bar\rho}(\epi f, \epi f^\nu) + \rho \|y-y^\nu\|_2, 
\]
where the norm on $\reals^m\times\reals^n \times \reals$ is $(u,x,\alpha)\mapsto \max\{\|u\|_2, \|x\|_2,|\alpha|\}$.
\end{enumerate}
\end{theorem}

If $f$ is a Rockafellian for $\phi:\reals^n\to \Reals$, then Theorem \ref{thm:errorRock}(a)  allows us to conclude that $\nlimsup (\inf f_{y^\nu}^\nu ) \leq \inf f_y \leq \inf \phi$; cf. Proposition \ref{prop:tigthepi}(b) below. Thus, the approximating Rockafellian relaxation is indeed a relaxation of the actual problem in the limit regardless of $y$ as long as $y^\nu\to y$. The rate of convergence of minimum values, near-minimizers, and level-sets of approximating tilted Rockafellians to those of the tilted Rockafellian follows from Theorem \ref{thm:errorRock}(b) and Proposition \ref{tapproxoptimalvalue}. These facts hold broadly in sharp contrast to the situation for $\phi^\nu$ relative to $\phi$ and thus justify the consideration of Rockafellian relaxations.

\begin{figure}
\centering
\includegraphics[width=0.7\textwidth]{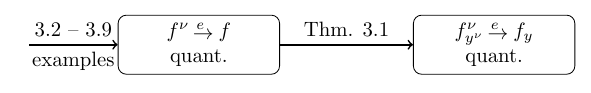}
\caption{Usage of Theorem \ref{thm:errorRock} to establish $f^\nu_{y^\nu}\eto f_y$ and its quantification.}\label{overview3}
\end{figure}

In the remainder of the section, we show how $f^\nu \eto f$ and its quantification by the truncated Hausdorff distance may hold in various applications. This provides the necessary input to Theorem \ref{thm:errorRock}; see Figure \ref{overview3} for an overview.  

\subsection{Epi-Convergence and Quantification in Composite Optimization}\label{subsec:compositeRock}

Theorem \ref{thm:errorRock} assumes epi-convergence of Rockafellians, which holds naturally in the setting of Example \ref{e:composite}. A short proof appears in Section \ref{sec:supp}.  

\begin{proposition}{\rm (composite optimization; epi-convergence).}\label{prop:composite}
In the setting of Example \ref{e:composite}, suppose that $X^\nu\sto X$, $h^\nu\eto h$, with $h$ being proper, and $g_0^\nu(x^\nu)\to g_0(x)$ and $G^\nu(x^\nu)\to G(x)$ whenever $x^\nu\in X^\nu\to x\in X$. Then, $f^\nu\eto f$.
\end{proposition}

The assumptions of Proposition \ref{prop:composite} fail to ensure that the composite functions $\phi^\nu$ epi-converge to $\phi$. The simple nonlinear program in Example \ref{e:composite} furnishes a counterexample. Still, the proposition confirms the epi-convergence of the Rockafellians of Example \ref{e:composite} and thus also of the corresponding tilted Rockafellians by Theorem \ref{thm:errorRock}(a). Informally, this means that small changes to components of a composite problem have only small effect on near-minimizers, minima, and level-sets of Rockafellian relaxations. In fact, the effect can be quantified in a nonasymptotic analysis as seen next for a fixed $\nu$.

\begin{proposition}{\rm (composite optimization; truncated Hausdorff distance).}\label{prop:compositeHaus}
In the setting of Example \ref{e:composite}, suppose that $X, X^\nu$ are nonempty, $h, h^\nu$ are proper, and $g_i, g_i^\nu:\reals^n\to \reals$, $i=0, 1, \dots, m$, are Lipschitz with common modulus function $\kappa$, where $g_1, \dots, g_m$ are the component functions of $G$ and $g_1^\nu, \dots, g_m^\nu$ are the component functions of $G^\nu$. Then, for any $\rho \in [0,\infty)$, one has
\begin{align*}
\hatsetd_\rho(\epi f, \epi f^\nu) \leq & \max\{1, \sqrt{m}\kappa(\rho')\} \hatsetd_\rho(X, X^\nu) + \hatsetd_{\bar\rho}(\epi h, \epi h^\nu)\\
&~~ + \sqrt{m}\max_{i=0, 1, \dots, m} \sup_{x\in \ball(\rho)} \big|g_i(x) - g_i^\nu(x)\big|
\end{align*}
provided that $\rho' > \rho + \hatsetd_\rho(X,X^\nu)$ and 
\[
\bar\rho \geq \rho + \sup_{x\in \ball(\rho)} \max\big\{|g_0(x)|, |g_0^\nu(x)|, \|G(x)\|_2, \|G^\nu(x)\|_2\big\}. 
\]
Here, the norms on $\reals^n$, $\reals^m \times \reals$, and $\reals^m\times\reals^n\times\reals$ are $\|\cdot\|_2$, $(z,\alpha)\mapsto \max\{\|z\|_2, |\alpha|\}$, and $(u,x,\alpha)\mapsto \max\{\|u\|_2, \|x\|_2, |\alpha|\}$, respectively. 
\end{proposition}
\state Proof. We adopt the following notation 
\begin{align*}
\delta_0 &= \sup_{x\in \ball(\rho)} \max\big\{|g_0(x)|, |g_0^\nu(x)|\big\}, ~~~~\delta = \sup_{x\in \ball(\rho)} \max\big\{\|G(x)\|_2, \|G^\nu(x)\|_2\big\}\\
\eta &=  \sup_{x\in \ball(\rho)} \max_{i=0, \dots, m}\big|g_i(x) - g_i^\nu(x)\big|.
\end{align*}
Let $\epsilon \in (0, \rho' - \rho - \hatsetd_\rho(X,X^\nu)]$ and $(u^\nu,x^\nu,\alpha^\nu) \in \epi f^\nu \cap \ball(\rho)$. Thus, $x^\nu\in X^\nu$ because $f^\nu(u^\nu,x^\nu) \leq \alpha^\nu\leq \rho$. Since $\|x^\nu\|_2 \leq \rho$ as well, there exists $x\in X$ such that $\|x - x^\nu\|_2 \leq \hatsetd_\rho(X,X^\nu) + \epsilon$. Moreover, $f^\nu(u^\nu,x^\nu) \leq \alpha^\nu$ implies that
$h^\nu\big(G^\nu(x^\nu) + u^\nu\big) \leq \alpha^\nu - g_0^\nu(x^\nu)$. Since $\|G^\nu(x^\nu) + u^\nu\|_2 \leq \delta + \rho$ and $|\alpha^\nu - g_0^\nu(x^\nu)| \leq \rho + \delta_0$, one has $(G^\nu(x^\nu) + u^\nu, ~\alpha^\nu - g_0^\nu(x^\nu)) \in \epi h^\nu \cap \ball(\bar\rho)$. By the definition of the truncated Hausdorff distance, we can find $(z,\beta)\in \epi h$ such that  
\[
\max\big\{\|G^\nu(x^\nu) + u^\nu - z\|_2, |\alpha^\nu - g_0^\nu(x^\nu) - \beta|\big\} \leq \hatsetd_{\bar\rho}(\epi h, \epi h^\nu) + \epsilon.
\] 
We set $u = z - G(x)$, which yields
\begin{align*}
  f(u,x) & = g_0(x) + h(z) \leq g_0(x) + \beta  \leq g_0(x) + \alpha^\nu - g_0^\nu(x^\nu) + \hatsetd_{\bar\rho}(\epi h, \epi h^\nu) + \epsilon\\
         & \leq \big|g_0(x) - g_0(x^\nu)\big| + \big|g_0(x^\nu) - g_0^\nu(x^\nu)\big| + \alpha^\nu + \hatsetd_{\bar\rho}(\epi h, \epi h^\nu) + \epsilon\\
         & \leq \kappa(\rho')\hatsetd_\rho(X,X^\nu) + \kappa(\rho')\epsilon + \eta + \alpha^\nu + \hatsetd_{\bar\rho}(\epi h, \epi h^\nu) + \epsilon
\end{align*}
because $\|x\|_2 \leq \|x^\nu\|_2 + \|x-x^\nu\|_2 \leq \rho + \hatsetd_\rho(X,X^\nu) + \epsilon \leq \rho'$. Let 
\[
\alpha = \kappa(\rho')\hatsetd_\rho(X,X^\nu) + \kappa(\rho')\epsilon + \eta + \alpha^\nu + \hatsetd_{\bar\rho}(\epi h, \epi h^\nu) + \epsilon.
\]
We also obtain the bounds
\begin{align*}
  \|u - u^\nu\|_2 & \leq \big\|z - G(x) - u^\nu\big\|_2\\
   & \leq \big\|G(x) - G(x^\nu)\big\|_2 + \big\|G(x^\nu) - G^\nu(x^\nu)\big\|_2 + \big\|G^\nu(x^\nu) + u^\nu - z\big\|_2\\
   & \leq \sqrt{m}\kappa(\rho')\hatsetd_\rho(X,X^\nu) + \sqrt{m}\kappa(\rho')\epsilon + \sqrt{m}\eta +   \hatsetd_{\bar\rho}(\epi h, \epi h^\nu) + \epsilon;\\
 |\alpha - \alpha^\nu| & \leq  \kappa(\rho')\hatsetd_\rho(X,X^\nu) + \kappa(\rho')\epsilon + \eta + \hatsetd_{\bar\rho}(\epi h, \epi h^\nu) + \epsilon.
\end{align*}
We have constructed $(u,x,\alpha) \in \epi f$ with distance to $(u^\nu,x^\nu,\alpha^\nu)$ of at most
\begin{align*}
\max\Big\{&\sqrt{m}\kappa(\rho')\hatsetd_\rho(X,X^\nu) + \sqrt{m}\kappa(\rho')\epsilon + \sqrt{m}\eta +   \hatsetd_{\bar\rho}(\epi h, \epi h^\nu) + \epsilon,\\
& \hatsetd_\rho(X,X^\nu) + \epsilon, ~~ \kappa(\rho')\hatsetd_\rho(X,X^\nu) + \kappa(\rho')\epsilon + \eta + \hatsetd_{\bar\rho}(\epi h, \epi h^\nu) + \epsilon\Big\}.
\end{align*}
Since $(u^\nu,x^\nu,\alpha^\nu)$ is arbitrary, this bound is also valid for $\exs(\epi f^\nu \cap \ball(\rho); \epi f)$. We repeat this argument with the roles of $f$ and $f^\nu$ reversed. Thus, $\hatsetd_\rho(\epi f, \epi f^\nu)$ is bounded by the same quantity. By letting $\epsilon$ tend to zero and slightly simplify the expression, we reach the conclusion.\eop

The proposition furnishes the key input to Theorem \ref{thm:errorRock}(b). Together, they bound $\hatsetd_\rho(\epi f_y, \epi f^\nu_{y^\nu})$ in terms of the discrepancies between the respective components of the composite problems. In turn, the bound feeds into Proposition \ref{tapproxoptimalvalue}. This takes place in a setting where there is no guarantee that   $\hatsetd_\rho(\epi \phi, \epi \phi^\nu)$ is small. 

The Rockafellian $f$ in Example \ref{e:composite} may not be exact but Subsection \ref{subsec:augmentationRock} provides a remedy.

\subsection{Epi-Convergence and Quantification under Constraint Structure}\label{subsec:constraint}

Beyond Example \ref{e:composite}, we also consider two examples with constraint structure to further illustrate the presence of epi-convergence and the calculation of $\hatsetd_\rho(\epi f, \epi f^\nu)$, which again allow for the use of Theorem \ref{thm:errorRock}.

\begin{example}{\rm (inequalities).}\label{e:constrStucture} For proper functions $\{g_i,g_i^\nu:\reals^n\to \Reals$, $i = 0, 1$, $\dots$, $m$, $\nu\in\nats\}$, consider
\[
\phi(x) = g_0(x) + \nsum_{i=1}^m \iota_{(-\infty,0]}\big(g_i(x)\big) ~\mbox{ and }~ \phi^\nu(x) = g_0^\nu(x) + \nsum_{i=1}^m \iota_{(-\infty,0]}\big(g_i^\nu(x)\big).
\]
A perturbation vector $u =(v_0, v_1, \dots, v_m, w_1, \dots, w_m)$, with $v_i \in \reals^n$ and $w_i \in \reals$, defines Rockafellians $f,f^\nu:\reals^{n(m+1) + m} \times \reals^n \to \Reals$ for $\phi$ and $\phi^\nu$, respectively, as 
\begin{align*}
f(u,x) & = g_0(x+v_0) + \nsum_{i=1}^m \iota_{(-\infty, 0]}\big(g_i(x+v_i) + w_i\big)\\
f^\nu(u,x) & = g_0^\nu(x+v_0) + \nsum_{i=1}^m \iota_{(-\infty, 0]}\big(g_i^\nu(x+v_i) + w_i\big).
\end{align*}
If $g_i^\nu \eto g_i$, $i = 0, 1, \dots, m$, then $f^\nu \eto f$. Regardless of the epi-convergence, for any $\nu\in \nats$ and $\rho \in [0,\infty)$, 
\[
\hatsetd_{\rho}(\epi f, \epi f^\nu) \leq \max_{i=0, 1,\dots, m} \hatsetd_{2\rho}(\epi g_i, \epi g_i^\nu),
\]
where the norms on $\reals^n\times \reals$ and $\reals^{n(m+1) + m} \times \reals^n \times \reals$ are $(x,\alpha)\mapsto \max\{\|x\|_2, |\alpha|\}$ and 
\[
(u,x,\alpha)\mapsto \max\big\{\|v_0\|_2, \|v_1\|_2, \dots, \|v_m\|_2, |w_1|, \dots, |w_m|, \|x\|_2, |\alpha|\big\},
\]
respectively, for $u =(v_0, v_1, \dots, v_m, w_1, \dots, w_m)$. 

In the simple nonlinear program of Example \ref{e:composite}, $\hatsetd_{2\rho}(\epi g_0, \epi g_0^\nu) = 0$ and $\hatsetd_{2\rho}(\epi g_1, \epi g_1^\nu)$ $=$ $1/\nu$ so that $\hatsetd_{\rho}(\epi f, \epi f^\nu) \leq 1/\nu$. In contrast, $\hatsetd_\rho(\epi \phi,\epi \phi^\nu)$ remains large for all $\nu$. 
\end{example}
\state Detail. The liminf-condition in \eqref{eqn:liminf} follows immediately. For the limsup-condition \eqref{eqn:limsup}, let $x\in\reals^n$, $v_i\in \reals^n$, $i=0, 1, \dots, m$, and $w\in \reals^m$. Without loss of generality, we assume that $g_i(x + v_i) + w_i \leq 0$ for $i=1, \dots, m$. For each $i$, there exists $z_i^\nu \to x + v_i$ such that $g_i^\nu(z_i^\nu) \to g_i(x+v_i)$. Construct $x^\nu = x$, $v_i^\nu = z_i^\nu - x$, $i=0, 1, \dots, m$, and $w_i^\nu = w_i - \max\{0, g_i^\nu(x^\nu + v_i^\nu) - g_i(x + v_i)\}$, $i = 1,\dots, m$. Then, 
\begin{align*}
 & \nlimsup f^\nu\big((v_0^\nu, \dots, v_m^\nu,w_1^\nu, \dots, w_m^\nu),x^\nu\big)\\
  & \leq \nlimsup g_0^\nu(z_0^\nu) + \nsum_{i=1}^m \nlimsup \iota_{(-\infty, 0]}\big(g_i^\nu(z_i^\nu) + w_i^\nu\big) = g_0(x + v_0). 
\end{align*}
This confirms the claim about epi-convergence.

For the quantification, let $\epsilon \in (0,\infty)$ and $(u,x,\alpha) \in \epi f \cap \ball(\rho)$, where $u =(v_0, v_1, \dots, v_m, w_1, \dots, w_m)$. Since $(x+v_0,\alpha)\in \epi g_0$, $\|x + v_0\|_2 \leq \|x\|_2 + \|v_0\|_2 \leq 2\rho$, $|\alpha|\leq \rho$, and $\epi g_0^\nu$ is nonempty, there exist $\bar v_0\in \reals^n$ and $\bar \alpha\in \reals$ such that $(x + \bar v_0, \bar\alpha) \in \epi g_0^\nu$ and
\[
\max\big\{\|x + \bar v_0 - x - v_0\|_2, |\bar\alpha - \alpha|\big\} \leq \hatsetd_{2\rho}(\epi g_0, \epi g_0^\nu) + \epsilon.
\]
Moreover, for $i = 1, \dots, m$, we have that $(x+v_i,-w_i) \in \epi g_i$, $\|x+v_i\|_2 \leq 2\rho$, and $|w_i| \leq \rho$. Since $\epi g_i^\nu$ is nonempty, there exist $\bar v_i\in \reals^n$ and $\hat \alpha_i\in \reals$ such that $(x + \bar v_i, \hat\alpha_i) \in \epi g_i^\nu$ and
\[
\max\big\{\|x + \bar v_i - x - v_i\|_2, |\hat\alpha_i + w_i|\big\} \leq \hatsetd_{2\rho}(\epi g_i, \epi g_i^\nu) + \epsilon.
\]
We construct $\bar x = x$ and $\bar u = (\bar v_0, \bar v_1, \dots, \bar v_m, \bar w_1, \dots, \bar w_m)$, where $\bar w_i = -\hat\alpha_i$, $i=1, \dots, m$. Then, $g_i^\nu(\bar x + \bar v_i) + \bar w_i = g_i^\nu(x + \bar v_i) - \hat \alpha_i \leq 0$ for $i=1, \dots, m$. This implies that $f^\nu(\bar u, \bar x) = g_0^\nu(x + \bar v_0) \leq \bar \alpha$. Consequently, $(\bar u,  \bar x, \bar\alpha) \in \epi f^\nu$. The distance between $(u, x, \alpha)$ and $(\bar u,  \bar x, \bar\alpha)$ is determined by
\begin{align*}
    &\|\bar v_i - v_i\|_2 \leq  \hatsetd_{2\rho}(\epi g_i, \epi g_i^\nu) + \epsilon, ~~i = 0, 1, \dots, m\\
    &|\bar w_i - w_i| = |-\hat\alpha_i - w_i| \leq \hatsetd_{2\rho}(\epi g_i, \epi g_i^\nu) + \epsilon, ~~i=1, \dots, m\\
    &|\bar \alpha - \alpha| \leq \hatsetd_{2\rho}(\epi g_0, \epi g_0^\nu) + \epsilon.
\end{align*}
The choice of norms and the fact that $\epsilon$ is arbitrary lead to the conclusion.\eop

\begin{example}{\rm (composite structure in constraint).}\label{e:compositeInConstr}
For proper $g_0,g_0^\nu:\reals^n\to \Reals$, proper $h,h^\nu:\reals^m\to \Reals$, and $G,G^\nu:\reals^n\to \reals^m$, $\nu\in\nats$, consider
\[
\phi(x) = g_0(x) + \iota_{(-\infty,0]}\big( h\big(G(x)\big)\big) ~~\mbox{ and }~~ \phi^\nu(x) = g_0^\nu(x) + \iota_{(-\infty,0]}\big( h^\nu\big(G^\nu(x)\big)\big).
\]
A perturbation vector $u = (v,w)$, with $v \in \reals^m$ and $w \in \reals$, defines Rockafellians $f,f^\nu:\reals^{m+1} \times \reals^n \to \Reals$ for $\phi$ and $\phi^\nu$, respectively, as 
\begin{align*}
f(u,x) & = g_0(x) + \iota_{(-\infty,0]}\big( h\big(G(x) + v\big) + w \big)\\
f^\nu(u,x) & = g_0^\nu(x) + \iota_{(-\infty,0]}\big( h^\nu\big(G^\nu(x) + v\big) + w\big). 
\end{align*}
Then, arguments similar to those in the proof of Proposition \ref{prop:compositeHaus} lead to:
\begin{enumerate}[(a)]

\item If  $g^\nu_0 \eto g_0$, $h^\nu \eto h$, and $G^\nu(x^\nu) \to G(x)$ for all $x^\nu\in \dom g_0^\nu\to x\in \dom g_0$, then $f^\nu\eto f$. 

\item If $\rho \in [0,\infty)$ and each component function $g_i$ and $g_i^\nu$ of $G$ and $G^\nu$, respectively, are Lipschitz with common modulus function $\kappa$, then
\begin{align*}
\hatsetd_\rho(\epi f^\nu, \epi f) \leq & \max\big\{1,\sqrt{m}\kappa(\rho')\big\}\hatsetd_{\rho}(\epi g_0, \epi g_0^\nu) + \hatsetd_{\bar\rho}(\epi h, \epi h^\nu)\\
&~~ + \sup_{\|x\|_2\leq \rho'}\big\|G(x) - G^\nu(x)\big\|_2
\end{align*}
provided that $\rho' > \rho + \hatsetd_\rho(\epi g_0, \epi g_0^\nu)$ and 
\[
\bar\rho \geq \rho + \sup_{\|x\|_2\leq \rho} \max\big\{\|G(x)\|_2,\|G^\nu(x)\|_2\big\}.
\]  
Here, we adopt the norms $(x,\alpha)\mapsto \max\{\|x\|_2,|\alpha|\}$ on $\reals^n \times \reals$,  $(z,\alpha)\mapsto \max\{\|z\|_2,$ $|\alpha|\}$ on $\reals^m\times\reals$, and $(u,x,\alpha) \mapsto \max\{\|v\|_2, |w|, \|x\|_2, |\alpha|\}$ on $\reals^{m+1}\times\reals^n\times\reals$, where $u = (v,w)$.
\end{enumerate}
\end{example}

\subsection{\redrev{Epi-Convergence and Quantification for} Stochastic Optimization}

Data-driven stochastic optimization gives rise to approximations stemming from imprecise probability models. 
Let $\Delta = \{p\in [0,\infty)^m~|~\nsum_{i=1}^m p_i = 1\}$ be the $m$-dimensional probability simplex, with $p\in \Delta$ specifying the weight assigned to each data point. Small changes to $p$ may have a dramatic effect on the solutions of the corresponding problems \cite{RoysetChenEckstrand.22}, but this can be mitigated by considering Rockafellian relaxations. In fact, $\phi^\nu$ may not epi-converge to $\phi$ in the next propositions but novel Rockafellians do.

\begin{proposition}{\rm (ambiguity).}\label{e:ambiguity}
For proper lsc $\{g_i:\reals^n\to \Reals, i=0, 1, \dots, m\}$ and $\{p,p^\nu\in \Delta, \nu\in\nats\}$, let 
\[
\phi(x) = g_0(x) + \nsum_{i=1}^m p_i g_i(x)~~~~\mbox{ and }~~~~\phi^\nu(x) = g_0(x) + \nsum_{i=1}^m p_i^\nu g_i(x).
\]
We adopt the Rockafellians $f,f^\nu:\reals^m\times\reals^n\to \Reals$ for $\phi$ and $\phi^\nu$, respectively, given by
\begin{align*}
f(u,x) &= g_0(x) + \nsum_{i=1}^m (p_i + u_i) g_i(x) + \half\theta \|u\|_2^2 + \iota_{(-\infty,0]^m}(u) + \iota_{[0,\infty)^m}(p+u)\\
f^\nu(u,x) &= g_0(x) + \nsum_{i=1}^m (p_i^\nu + u_i) g_i(x) + \half\theta^\nu \|u\|_2^2 + \iota_{(-\infty,0]^m}(u) + \iota_{[0,\infty)^m}(p^\nu+u),
\end{align*}
where $\theta, \theta^\nu \in [0,\infty)$. Then, the following hold: 
\begin{enumerate}[(a)]

\item If $\ninf g_i \in \reals$\redrev{, i = 0, 1, \ldots, m}, and there is $\bar x\in\reals^n$ such that $g_i(\bar x) <\infty$ for $i = 0, 1, \dots, m$, then there is $\bar y\in \reals^m$ such that, regardless of $\theta \in [0,\infty)$, the Rockafellian $f$ is strictly exact, supported by any $y\geq \bar y$.

\item If $p^\nu\to p$ and $\theta^\nu\to \theta$, then $f^\nu\eto f$.

\item If $\rho \in [0,\infty)$ and there is $\eta\in [0,\infty)$ such that $\ninf g_i \geq -\eta$ for $i=1, \dots, m$, then  
\[
\hatsetd_\rho(\epi f, \epi f^\nu) \leq \max\big\{1, \eta + \theta^\nu\}\|p-p^\nu\|_1 + \tfrac{1}{2}\theta^\nu\|p-p^\nu\|_2^2 + \tfrac{1}{2}|\theta - \theta^\nu|
\]
under the norm $(u,x,\alpha)\mapsto \max\{\|u\|_1, \|x\|_2,|\alpha|\}$ on $\reals^m\times\reals^n\times \reals$.
\end{enumerate}
\end{proposition}
\state Proof. Consider (a). For $u\in \reals^m$, we adopt the notation
\[
\delta(u) = \ninf_x g_0(x) + \nsum_{i=1}^m (p_i + u_i) g_i(x) + \iota_{(-\infty,0]^m}(u) + \iota_{[0,\infty)^m}(p+u).
\]
By assumption, there is $\bar \eta \in [0,\infty)$ such that $\ninf g_i \geq -\bar \eta$ for all $i = 0, 1, \dots, m$. Certainly,
\[
-2\bar \eta \leq \delta(0) \leq g_0(\bar x) + \nsum_{i=1}^m p_i g_i(\bar x)<\infty.
\]
We adopt the notation
$I_0 = \{i~|~p_i =0\}$, $I^+ = \{i~|~p_i > 0\}$, and $\alpha = \min\{p_i~|~i\in I^+\}$.

Suppose that $u = (u_1, \dots, u_m)$ satisfies $\max\{-\alpha/2, - p_i\} \leq u_i \leq 0$. We consider two cases: 
Case A: Suppose that $\delta(u) < \delta(0)$. Let $\epsilon \in (0,\delta(0) - \delta(u)]$. \redrev{Then, $\delta(u) \geq -2\bar\eta$ because} 
\[
\redrev{\delta(u) \geq -\bar\eta - \nsum_{i=1}^m (p_i + u_i) \bar\eta \geq -2\bar\eta.}
\]
Since $-2\bar \eta \leq \delta(u) < \delta(0)<\infty$ \redrev{and the indicator functions in the definition of $f$ vanish at $u$}, there exists $\hat x\in\reals^n$ such that
\[
g_0(\hat x) + \nsum_{i=1}^m (p_i + u_i) g_i(\hat x) \leq \delta(u) + \epsilon \leq \delta(0).
\]
Let $i\in I^+$. Then, 
\[
(p_i + u_i) g_i(\hat x) \leq \delta(0) - g_0(\hat x) - \nsum_{j\neq i} (p_j + u_j) g_j(\hat x) \leq \delta(0) + 2\bar \eta.
\] 
The sum $p_i + u_i\geq \alpha/2$ because $p_i\geq \alpha$ and $u_i \in [-\alpha/2, 0]$. Consequently, $g_i(\hat x) \leq (\delta(0) + 2\bar \eta)/(p_i + u_i) \leq (2\delta(0) + 4\bar \eta)/\alpha<\infty$. This bound together with $g_i(\hat x)\geq -\bar \eta$ \redrev{and the fact that $u_i = 0$ for $i\in I_0$} facilitate the following development:
\begin{align*}
  \delta(u) &\geq g_0(\hat x) + \nsum_{i \in I^+} (p_i + u_i) g_i(\hat x) + \nsum_{i \in I_0} (p_i + u_i) g_i(\hat x) - \epsilon\\
   & \geq g_0(\hat x) + \nsum_{i \in I^+} p_i g_i(\hat x) + \nsum_{i \in I^+}  u_i g_i(\hat x) - \epsilon\\
&\geq g_0(\hat x) + \nsum_{i \in I^+} p_i g_i(\hat x) + \nsum_{i \in I^+} u_i \frac{2\delta(0) + 4\bar \eta}{\alpha} - \epsilon\\
&\geq \ninf_x g_0(x) + \nsum_{i=1}^m p_i g_i(x) + \langle y, u\rangle - \epsilon,
\end{align*}
where $y = (y_1, \dots, y_m)$, with $y_i \geq (2\delta(0) + 4\bar \eta)/\alpha$ if $i \in I^+$ and $y_i \geq 0$ if $i\in I_0$. Since $\epsilon$ is arbitrary, we have shown that $\delta(u) \geq \delta(0) + \langle y, u\rangle$ when $\delta(u) < \delta(0)$. 

Case B: Suppose that $\delta(u) \geq \delta(0)$. Since $u\leq 0$ and $y\geq 0$, we trivially have that $\delta(u) \geq \delta(0) + \langle y, u\rangle$.

We have shown that $\delta(u) \geq \delta(0) + \langle y, u\rangle$ whenever $u = (u_1, \dots, u_m)$ satisfies $\max\{-\alpha/2, - p_i\}$ $\leq u_i$ $\leq 0$. 

For $u' = (u'_1, \dots, u'_m)$ with $-p_i \leq u'_i \leq 0$ for all $i$ suppose that there exists $i^\star$ such that $u'_{i^\star} < -\alpha/2$. Then,  $\delta(u') \geq - 2\bar \eta -\delta(0) - \langle y, u'\rangle + \delta(0) + \langle y, u'\rangle \geq \delta(0) + \langle y, u'\rangle$ because 
\begin{align*}
 - 2\bar \eta -\delta(0) - \langle y, u'\rangle & \geq  - 2\bar \eta -\delta(0) - \nsum_{i \in I^+} \frac{2\delta(0) + 4\bar \eta}{\alpha}u_i'\\
  & \geq - 2\bar \eta -\delta(0) - \frac{2\delta(0) + 4\bar \eta}{\alpha}\frac{-\alpha}{2} = 0.
\end{align*}
Thus, $\delta(u) \geq \delta(0) + \langle y, u\rangle$ for all $u\in \reals^m$ \redrev{because the only $u$-vectors not considered have $\delta(u) = \infty$}. A slight increase in $y_i$ for $i \in I^+$ ensures strict exactness.

For (b), first let $x^\nu\to x$ and $u^\nu\to u$. Then, $\nliminf g_0(x^\nu)\geq g_0(x)>-\infty$, $\nliminf \iota_{(-\infty,0]^m}(u^\nu)$ $\geq$ $\iota_{(-\infty,0]^m}(u)$, and $\nliminf \iota_{[0,\infty)^m}(p^\nu+u^\nu) \geq \iota_{[0,\infty)^m}(p+u)$ because $(-\infty,0]^m$ and $[0,\infty)^m$ are closed sets. Since $f^\nu(u^\nu,x^\nu) = \infty$ if $p^\nu + u^\nu \not\in [0,\infty)^m$, we assume without loss of generality that  $p^\nu + u^\nu \in [0,\infty)^m$ for all $\nu\in\nats$. This in turn implies that
\[
\nliminf \nsum_{i=1}^m (p_i^\nu + u_i^\nu) g_i(x^\nu) \geq \nsum_{i=1}^m \nliminf  (p_i^\nu + u_i^\nu) g_i(x^\nu) \geq \nsum_{i=1}^m (p_i + u_i) g_i(x).
\]
Combining these facts, we find that $\nliminf f^\nu(u^\nu,x^\nu) \geq f(u,x)$.

Second, let $(u,x)\in \reals^m \times \reals^n$. Without loss of generality, we assume that $u \in (-\infty, 0]^m$ and $p+u \in [0,\infty)^m$.  We construct $x^\nu = x$ and $u_i^\nu = \min\{0, p_i + u_i - p_i^\nu\}$. Thus, $u^\nu_i \leq 0$,  $u_i^\nu\to u_i$, and $p_i^\nu + u_i^\nu \geq 0$. The latter holds because if $p_i + u_i - p^\nu_i \geq 0$, then $u^\nu_i = 0$ and $p^\nu_i + u^\nu_i = p^\nu_i \geq 0$. If $p_i+u_i - p^\nu_i < 0$, then $u^\nu_i = p_i + u_i - p^\nu_i$ and $p^\nu_i + u^\nu_i = p_i + u_i \geq 0$. These facts yield
\begin{align*}
\nlimsup f^\nu(u^\nu,x^\nu) & \leq g_0(x) + \nsum_{i=1}^m \nlimsup \big(p_i^\nu + \min\{0, p_i + u_i - p_i^\nu\}\big) g_i(x)\\
& ~~+ \nlimsup \half \theta^\nu\nsum_{i=1}^m \big(\min\{0, p_i + u_i - p_i^\nu\}\big)^2 \leq f(u,x).
\end{align*}
Thus, we have confirmed both \eqref{eqn:liminf} and \eqref{eqn:limsup}.

For (c), let $(u,x,\alpha) \in \epi f \cap \ball(\rho)$. Then, $-p_i \leq u_i \leq 0$ for $i=1, \dots, m$ \redrev{because $f(u,x)\leq \alpha <\infty$}. Construct $\bar x = x$ and $\bar u_i = \min\{0, p_i - p_i^\nu + u_i\}$. Then,  $-p_i^\nu \leq \bar u_i \leq 0$ for $i=1, \dots, m$ and $\|u- \bar u\|_2 \leq \|p- p^\nu\|_2$. Moreover, 
\begin{align*}
    f^\nu(\bar u, \bar x) & = f(u,x) + \nsum_{i=1}^m \big(\min\{p_i^\nu, p_i + u_i\} - (p_i + u_i) \big) g_i(x)\\
     & ~~~~~~~+ \tfrac{1}{2}\nsum_{i=1}^m (\theta^\nu \min\{0, p_i - p_i^\nu + u_i\}^2 -  \theta u_i^2)\\
    & \leq \alpha + \eta \nsum_{i=1}^m|p_i^\nu - p_i| + \tfrac{1}{2}\nsum_{i=1}^m \big(\theta^\nu \min\{0, p_i - p_i^\nu + u_i\}^2  - \theta u_i^2\big)\\
    & \leq \alpha + \eta \|p^\nu - p\|_1 + \tfrac{1}{2}\nsum_{i=1}^m \big(\theta^\nu (p_i - p_i^\nu)^2 + 2\theta^\nu p_i |p_i - p_i^\nu| + |\theta^\nu - \theta| p_i^2\big).
\end{align*}
The first inequality is valid because when $p_i^\nu < p_i + u_i$, one has
\[
(p_i^\nu - p_i - u_i) g_i(x) = (p_i + u_i - p_i^\nu)\big(-g_i(x)\big) \leq (p_i + u_i - p_i^\nu) \eta \leq (p_i - p_i^\nu) \eta.
\]
The conclusion then follows \redrev{in the light of the symmetry between $f$ and $f^\nu$}.\eop

\redrev{The proof of part (a) is constructive and thus furnishes a formula for $\bar y$.}
Both this proposition and the next one consider novel Rockafellians compared to \cite{RoysetChenEckstrand.22}.

\begin{proposition}{\rm (splitting).}\label{e:splitting}
For proper $\{g_i,g_i^\nu:\reals^n\to \Reals, i=1, \dots, m, \nu\in\nats\}$ and $\{p,p^\nu\in \Delta, \nu\in\nats\}$, consider
\[
\phi(x) = \nsum_{i=1}^m p_i g_i(x) ~~\mbox{ and } ~~ \phi^\nu(x) = \nsum_{i=1}^m p_i^\nu g_i^\nu(x) 
\]
and, with $u = (u_1, \dots, u_m)\in \reals^{nm}$, the corresponding Rockafellians given by  
\[
f(u,x) = \nsum_{i=1}^m p_i g_i(x + u_i) ~~\mbox{ and } ~~~ f^\nu(u,x) = \nsum_{i=1}^m p_i^\nu g_i^\nu(x + u_i).
\]
Then, the following hold: 
\begin{enumerate}[(a)]

\item If $p^\nu\to p$ and, for each $i$, $g_i^\nu\eto g_i$ with $g_i$ being either real-valued or $p_i > 0$ or $p_i = p_i^\nu = 0$ for all $\nu$, then $f^\nu \eto f$. 

\item If $\rho\in [0,\infty)$ and there is $\eta \in [0,\infty)$ such that $\inf g_i \geq -\eta$ and $\inf g_i^\nu \geq -\eta$ for $i = 1, \dots, m$, then 
\[
\hatsetd_\rho(\epi f, \epi f^\nu) \leq \lambda \|p^\nu - p\|_2 + \max_{i = 1, \dots, m} \hatsetd_{\bar\rho} (\epi g_i, \epi g_i^\nu), 
\]
where we adopt the norms $(u,x,\alpha) \mapsto \max\{\|u_1\|_2, \dots, \|u_m\|_2, \|x\|_2, |\alpha|\}$ and $(x,\alpha)$ $\mapsto$ $\max\{\|x\|_2,|\alpha|\}$ on $\reals^{nm} \times \reals^n \times \reals$ and $\reals^n \times \reals$, respectively, and 
\begin{align*}
\lambda & = \sqrt{\nsum_{i=1}^m \big(\bar\rho + \hatsetd_{\bar\rho}(\epi g_i, \epi g_i^\nu)\big)^2}\\
\bar \rho &\geq \max\bigg\{2\rho, ~\frac{\rho + \eta}{\beta}, ~\sup_{\|x\|_2\leq 2\rho} \max\Big\{\max_{i \in I_0}  g_i(x),  \max_{i \in I^\nu}  g_i^\nu(x)\Big\}\bigg\}, 
\end{align*}
with $\beta = \min\{\min_{i\not\in I_0} p_i$, $\min_{i\not\in I^\nu} p_i^\nu\}$,  $I_0 = \{i \,|\, p_i = 0\}$, and $I^\nu = \{i \,|\, p_i^\nu = 0\}$.
\end{enumerate}

\end{proposition}
\state Proof. For (a), let $u^\nu=(u_1^\nu, \dots, u_m^\nu)\to u = (u_1, \dots, u_m)$ and $x^\nu\to x$. Since $\nliminf p_i^\nu g_i^\nu(x^\nu + u_i^\nu) \geq p_i g_i(x + u)$, the liminf-condition \eqref{eqn:liminf} follows. The limsup-condition \eqref{eqn:limsup} relies on the following argument. Let $(u,x)$ be arbitrary. For each $i$, there exists $z_i^\nu \to x + u_i$ such that  $g_i^\nu(z_i^\nu) \to g_i(x + u_i)$. Construct $x^\nu = x$ and $u_i^\nu = z_i^\nu - x$.  Under the stated assumptions, $\nlimsup p_i^\nu g_i^\nu(x^\nu + u_i^\nu) = \nlimsup p_i^\nu g_i^\nu(z_i^\nu)$ $=$ $ p_i g_i(x + u_i)$; \redrev{the potential difficulty of $p_i =0 $ and $g_i(x + u_i) =\infty$ is ruled out}. 

For (b), let $(u,x,\alpha) \in \epi f \cap \ball(\rho)$. For each $i$, one has $\|x+u_i\|_2 \leq 2\rho$ and 
\[
p_i g_i(x+u_i) \leq \rho - \nsum_{j\neq i} p_{j} g_{j}(x+u_{j}) \leq \rho + \eta.
\]
We next construct $(\bar u, \bar x, \bar \alpha) \in \epi f^\nu$. Set $\bar x = x$ and let $\epsilon \in (0,\infty)$. For each $i$, we consider two cases: Case A: Suppose that $i\not\in I_0$. Then, $g_i(x + u_i) \leq (\rho + \eta)/p_i\leq \bar\rho$ and $g_i(x+u_i) \geq - \eta \geq -\bar\rho$. Thus, $(x+u_i, g_i(x+u_i)) \in \epi g_i \cap \ball(\bar\rho)$.  Case B: Suppose that $i\in I_0$. Then, $|g_i(x + u_i)|\leq \bar\rho$ and, again, $(x+u_i, g_i(x+u_i)) \in \epi g_i \cap \ball(\bar\rho)$. In either case there exist $\bar u_i$ and $\bar\alpha_i$ such that $(x+ \bar u_i, \bar\alpha_i) \in \epi g_i^\nu$ and  
\[
\max\big\{\|u_i - \bar u_i\|_2, \big|g_i(x+u_i) - \bar \alpha_i\big|\big\} \leq \hatsetd_{\bar \rho}(\epi g_i, \epi g_i^\nu) + \epsilon.
\]
We then find that 
\begin{align*}
\nsum_{i=1}^m p_i^\nu g_i^\nu(\bar x + \bar u_i) & = \nsum_{i=1}^m p_i g_i(x + u_i) + \nsum_{i=1}^m (p_i^\nu - p_i) g_i^\nu(x + \bar u_i)\\
&~~~~~~ + \nsum_{i=1}^m p_i\big(g_i^\nu(x + \bar u_i) - g_i(x + u_i)\big)\\
& \leq \alpha + \|p^\nu - p\|_2 \sqrt{\nsum_{i=1}^m \big(\bar\rho + \hatsetd_{\bar\rho}(\epi g_i, \epi g_i^\nu) + \epsilon\big)^2}\\
&~~~~~~ + \nsum_{i=1}^m p_i \big(\hatsetd_{\bar\rho}(\epi g_i, \epi g_i^\nu) + \epsilon\big).
\end{align*}
Consequently, we can set $\bar \alpha$ to the quantity on the \redrev{right-hand side}. We repeat the argument with the roles of $f$ and $f^\nu$ reversed. The conclusion follows after letting $\epsilon$ tend to zero.\eop

\subsection{\redrev{Epi-Convergence and Quantification under} Augmentation}\label{subsec:augmentationRock}

A Rockafellian $f$ for $\phi$ that is not (strictly) exact can be augmented with a function of $u$ to produce another Rockafellian for $\phi$ with that property. One approach relies on norms to define augmenting functions, which relates closely to augmented Lagrangians \cite[Sec. 11.K]{VaAn}. We concentrate on the augmenting function  $\iota_{\{0\}^m}$. Trivially, the resulting Rockafellian $f + \iota_{\{0\}^m}$ for $\phi$ is strictly exact, supported by any $y\in \reals^m$. The computational challenges associated with this augmenting function are mitigated by considering an approximating function $a^\nu$, which may only introduce minor inaccuracies compared to those anticipated from approximating $f$ by $f^\nu$. 
This is formalized next, with a proof in Section \ref{sec:supp}.

\begin{theorem}{\rm (epi-convergence under augmentation).}\label{thm:augmentation}
For $\{a,a^\nu:\reals^m\to [0, \infty],$ $\nu\in\nats\}$ with $a(0) = a^\nu(0) = 0$ and proper $\{f,f^\nu:\reals^m\times\reals^n\to \Reals,$ $\nu\in\nats\}$, suppose that $f,f^\nu$ are Rockafellians for $\phi,\phi^\nu:\reals^n\to \Reals$, respectively, and let 
\[
\bar f(u,x) = f(u,x) + a(u) ~~\mbox{ and } ~~\bar f^\nu(u,x) = f^\nu(u,x) + a^\nu(u). 
\]
Then, $\bar f$ and $\bar f^\nu$ are Rockafellians for $\phi$ and $\phi^\nu$, respectively, and the following hold: 
\begin{enumerate}[(a)]
    \item If $f^\nu\eto f$, $a^\nu\eto a = \iota_{\{0\}^m}$, and for each $x\in \dom \phi$ there exist $u^\nu\to 0$ and $x^\nu\to x$ such that $f^\nu(u^\nu,x^\nu)\to f(0,x)$ and $a^\nu(u^\nu)\to 0$, then $\bar f^\nu\eto \bar f$. 

    \item If $f^\nu\eto f$ and $a^\nu(u^\nu) \to a(u)$ whenever $u^\nu\to u$, then $\bar f^\nu\eto \bar f$. 

    \item For $\rho\in [0,\infty)$, if there is $\eta \in [0,\infty)$ such that $\inf_{(u,x)\in \ball(\rho)} f(u,x) \geq -\eta$ and $\inf_{(u,x)\in \ball(\rho)} f^\nu(u,x) \geq -\eta$, and $a, a^\nu$ are real-valued Lipschitz functions with common modulus function $\kappa$, then 
\[
\hatsetd_{\rho}(\epi \bar f, \epi \bar f^\nu) \leq \big(1+\kappa(\rho')\big) \hatsetd_{\bar\rho}(\epi f, \epi f^\nu) + \nsup_{\|u\|_2\leq \rho'} \big|a^\nu(u) - a(u)\big|
\]
provided that $\bar\rho \geq \max\{\rho,\eta\}$ and $\rho' > \rho + \hatsetd_{\bar\rho}(\epi f, \epi f^\nu)$, where we utilize the norms $(u,x) \mapsto \max\{\|u\|_2, \|x\|_2\}$ and $(u,x,\alpha) \mapsto \max\{\|u\|_2, \|x\|_2, |\alpha|\}$ on $\reals^m\times\reals^n$ and $\reals^m\times\reals^n\times \reals$, respectively.

\end{enumerate}
\end{theorem}

While a real-valued $a$ leads to a quantifiable rate of convergence in Theorem \ref{thm:augmentation}(c), there appears to be no analogous rate result for $a = \iota_{\{0\}^m}$ and we settle for a mere convergence guarantee; see Theorem \ref{thm:augmentation}(a). The challenge emerges already in the case of inequality constraints as seen in \cite[Thm. 4.5]{Royset.20b}. 

If $a = \theta\|\cdot\|_2^2$, then $|a(\bar u) - a(u)| \leq 2\theta\rho'\|\bar u - u\|_2$ whenever $\|\bar u\|_2 \leq \rho'$ and $\|u\|_2 \leq \rho'$. Thus, with $a^\nu = \theta^\nu\|\cdot\|_2^2$, one can set $\kappa(\rho') = 2\max\{\theta,\theta^\nu\}\rho'$ in the theorem. 

Theorem \ref{thm:augmentation} supports the augmentation of the Rockafellians in Example \ref{e:composite}. Proposition \ref{prop:composite} and Theorem \ref{thm:augmentation}(b) immediately ensure epi-convergence of the augmented Rockafellians, for instance under augmentation function $\theta\|\cdot\|_2^2$. The rate of convergence follows by Proposition \ref{prop:compositeHaus} and Theorem \ref{thm:augmentation}(c). 

The assumptions of Theorem \ref{thm:augmentation}(a) hold widely as seen next; the proof appears in Section \ref{sec:supp}.

\begin{proposition}{\rm (composite optim. under augmentation; epi-convergence).}\label{prop:compositeAug}
In the setting of Example \ref{e:composite}, suppose that $X^\nu\sto X$, $h^\nu\eto h$ and are all proper functions, $h^\nu(u)\to h(u)$ for all $u\in \dom h$, the mapping $G$ is continuous, $g_0^\nu(x^\nu)\to g_0(x)$ whenever $x^\nu\in X^\nu \to x\in X$, and $C$ is an open set containing $X$ and $X^\nu$ for all $\nu$ sufficiently large. Given $\alpha \in (0,\infty)$, $\theta^\nu\in [0,\infty)$, any norm $\|\cdot\|$ on $\reals^m$, and the Rockafellians $f$ and $f^\nu$ from Example \ref{e:composite}, define augmented Rockafellians $\bar f$ and $\bar f^\nu$ by setting 
\[
\bar f(u,x) = f(u,x) + \iota_{\{0\}^m}(u) ~~\mbox{ and } ~~ \bar f^\nu(u,x) = f^\nu(u,x) + \theta^\nu\|u\|^\alpha.
\]
If $\theta^\nu\to \infty$ and $\theta^\nu \sup_{x\in C} \|G^\nu(x) - G(x)\|^\alpha\to 0$, then $\bar f^\nu \eto \bar f$.
\end{proposition}

\section{Convergence of Lagrangians}\label{sec:lagr}

We next focus on Lagrangians and the middle arrow in Figure \ref{overview}. Given that $f^\nu \eto f$ and $y^\nu\to y$, what additional properties would suffice for $l^\nu( \,\cdot\, , y^\nu) \eto l( \,\cdot\, ,y)$ to hold? In the convex setting, this is well understood even in reflexive Banach spaces and for the broader question of epi/hypo-convergence of $l^\nu$ to $l$ \cite{AtAW88:sad}; see also \cite{Attouch.86,AttouchAzeWets.86,Bagh.96}, with rates of convergence results in \cite{Aze.88}. (The general theory of epi/hypo-convergence as pioneered in \cite{AttouchWets.83} does not require convex-concave bifunctions.) In this section, we focus on the nonconvex setting and the ``primal'' variables, with the multipliers being viewed as parameters. We recall that the role of multipliers tends to be somewhat diminished in the nonconvex case when the underlying Rockafellian is often constructed via augmentation. This motivates the present focus, which brings forth {\em tightness} as the key property. 

\begin{theorem}{\rm (Lagrangian relaxations; epi-convergence).}\label{thm:epiLagrangians}
For Rockafellians $f$, $f^\nu:$ $\reals^m \times \reals^n\to \Reals$, suppose that $f^\nu \eto f$ and multiplier vectors $y^\nu\to y\in \reals^m$. If $\{f^\nu_{y^\nu}(\,\cdot\,,x^\nu), \nu\in\nats\}$ is tight for all convergent sequence $\{x^\nu\in\reals^n, \nu\in\nats\}$, then  $l^\nu(\,\cdot\,,y^\nu) \eto l(\,\cdot\,,y)$.
\end{theorem}
\state Proof. We first confirm the liminf-condition in \eqref{eqn:liminf}. Let $x^\nu\to x$. Since $f^\nu_{y^\nu} \eto f_y$ by Theorem \ref{thm:errorRock}(a), it follows from \eqref{eqn:liminf} that $\nliminf f^\nu_{y^\nu}(u^\nu,x^\nu) \geq f_y(u,x)$ whenever $u^\nu\to u$. We can then invoke Proposition \ref{prop:tigthepi}(a) to conclude that $\nliminf l^\nu(x^\nu,y^\nu) = \nliminf (\inf f^\nu_{y^\nu}(\,\cdot\,,x^\nu)) \geq  \inf f_y(\,\cdot\,,x) = l(x,y)$. 

Second, we confirm the limsup-condition in \eqref{eqn:limsup}. Let $x\in\reals^n$. We construct $x^\nu\to x$ such that
\begin{equation}\label{eqn:limsupcondproofLagr}
\nlimsup l^\nu(x^\nu,y^\nu) \leq l(x,y)
\end{equation}
and consider three cases. Let $\epsilon \in (0,\infty)$.  (i) Suppose that $l(x,y) \in \reals$. Then, there exists $u\in \epsilon\mbox{-}\nargmin f_y(\,\cdot\,,x)$. Since $f^\nu_{y^\nu}\eto f_y$, there are points $(u^\nu,x^\nu)\to (u,x)$ such that $\nlimsup f^\nu_{y^\nu}(u^\nu,x^\nu) \leq f_y(u,x)$. Thus,
\begin{align*}
\nlimsup l^\nu(x^\nu,y^\nu) & = \nlimsup \big(\inf f^\nu_{y^\nu}(\,\cdot\,,x^\nu)\big) \leq  \nlimsup f^\nu_{y^\nu}(u^\nu,x^\nu)\\
& \leq f_y(u,x) \leq \inf f_y(\,\cdot\,,x)+\epsilon = l(x,y)+\epsilon.
\end{align*}
Since $\epsilon$ is arbitrary, \eqref{eqn:limsupcondproofLagr} holds. (ii) The case with $l(x,y) = -\infty$ follows by a similar argument. (iii) Suppose that $l(x,y) = \infty$. Then, \eqref{eqn:limsupcondproofLagr} holds trivially with $x^\nu = x$ for all $\nu$.\eop

A consequence of the theorem is that $\nlimsup (\inf l^\nu(\,\cdot\,, y^\nu) )$ $\leq$ $\inf f(0, \,\cdot\,)$ $=$ $\inf \phi$; see Proposition \ref{prop:tigthepi}(b) below. Consequently, under the assumptions of Theorem \ref{thm:epiLagrangians} the approximating Lagrangians provide in the limit a lower bound on the minimum value of the actual problem.

We note that tightness is a necessary assumption in the theorem in the following sense. If there exists a convergent sequence $x^\nu\to x$ for which $\{f^\nu_{y^\nu}(\,\cdot\,,x^\nu), \nu\in\nats\}$ is not tight, then $\nliminf l^\nu(x^\nu,y^\nu) < l(x,y)$ and the epi-convergence of $l^\nu(\,\cdot\,,y^\nu)$ to $l(\,\cdot\,,y)$ will not hold. 

A result related to Theorem \ref{thm:epiLagrangians} appears as Exercise 7.57 in \cite{VaAn}. It asserts the following: For proper lsc $f,f^\nu$, with $f^\infty(u,0) > 0$ for all $u\neq 0$, we find that total epi-convergence of $f^\nu$ to $f$ implies that $l^\nu(\,\cdot\,,0)$ totally epi-converges to $l(\,\cdot\,,0)$. (Here, $f^\infty$ is the horizon function of $f$; see \cite[Sec. 3.C]{VaAn}.) We recall that total epi-convergence is a stronger property than epi-convergence but they are equivalent in some cases, for example for convex functions; see \cite[Thm. 7.53]{VaAn}. Our assumption about tightness in Theorem \ref{thm:epiLagrangians} is weaker than those of Exercise 7.57. Consider the example with $f(u,x) = f^\nu(u,x) = \max\{0,u\}$. Trivially, $f^\nu\eto f$, and this also holds in the sense of total epi-convergence because the functions are convex. For any convergent sequence $\{x^\nu, \nu\in\nats\}$,  $\{f^\nu(\,\cdot\,,x^\nu), \nu\in\nats\}$ is tight. Thus, Theorem \ref{thm:epiLagrangians} applies for $y^\nu = y = 0$. However, the requirement $f^\infty(u,0) > 0$ for all $u\neq 0$ in  \cite[Exercise 7.57]{VaAn} fails because the horizon function of $f^\infty = f$ in this case. 

We also establish a bound on $\hatsetd_\rho(\epi l(\,\cdot\,,y), \epi l^\nu(\,\cdot\,,y^\nu))$ that resembles the one in Theorem \ref{thm:errorRock}.

\begin{theorem}{\rm (Lagrangian relaxations; error bounds).}\label{thm:errorLagrangians}
For Rockafellians $f,f^\nu:\reals^m \times \reals^n\to \Reals$, suppose that $\epi f$ and $\epi f^\nu$ are nonempty, $\rho \in [0,\infty)$, and $y,y^\nu\in\reals^m$. If $\hat\rho\in [\rho, \infty)$ is such that
\begin{align*}
&\forall x\in \reals^n \mbox{ with } \|x\|_2\leq \rho, ~l(x,y) \leq \rho:\\
&~~~~\exists \hat u\in \reals^m \mbox{ with } \|\hat u\|_2 \leq \hat\rho, ~~f_y(\hat u,x) \leq \max\big\{-\rho, ~l(x,y)\big\},\\
&\forall x\in \reals^n \mbox{ with } \|x\|_2\leq \rho, ~l^\nu(x,y^\nu) \leq \rho\\
&~~~~\exists \hat u^\nu\in \reals^m \mbox{ with } \|\hat u^\nu\|_2 \leq \hat\rho, ~f^\nu_{y^\nu}(\hat u^\nu,x) \leq \max\big\{-\rho, ~l^\nu(x,y^\nu)\big\},
\end{align*}
then, using the norms $(u,x,\alpha)\mapsto \max\{\|u\|_2, \|x\|_2,|\alpha|\}$ and $(x,\alpha)\mapsto \max\{\|x\|_2,|\alpha|\}$ on  $\reals^m\times\reals^n \times \reals$ and $\reals^n \times \reals$, respectively, and with $\rho' \geq (1 + \max\{\|y\|_2, \|y^\nu\|_2\})\hat\rho$, one has
\[
\hatsetd_\rho\big(\epi l(\,\cdot\,,y), \epi l^\nu(\,\cdot\,,y^\nu)\big) \leq \Big( 1 + \max\big\{\|y\|_2, \|y^\nu\|_2\big\} \Big)\hatsetd_{\rho'}(\epi f, \epi f^\nu) + \hat\rho \|y-y^\nu\|_2.
\]
\end{theorem}
\state Proof. Set $\epsilon \in (0,\infty)$. Let $(x^\nu,\alpha^\nu)\in \epi l^\nu(\,\cdot\,, y^\nu) \cap \ball(\rho)$. Since $\|x^\nu\|_2 \leq \rho$ and $l^\nu(x^\nu,y^\nu) \leq \alpha^\nu\leq \rho$, there is $\hat u^\nu\in \ball(\hat\rho)$ with $f^\nu_{y^\nu}(\hat u^\nu,x^\nu) \leq \max\{-\rho, ~l^\nu(x^\nu,y^\nu)\}$. 

We consider two cases. First, suppose that $l^\nu(x^\nu,y^\nu) \geq -\rho$. Then, $f^\nu_{y^\nu}(\hat u^\nu,x^\nu) = l^\nu(x^\nu,y^\nu) \leq \alpha^\nu$. This means that $(\hat u^\nu,x^\nu,\alpha^\nu)$ $\in$ $\epi f^\nu_{y^\nu} \cap \ball(\hat\rho)$.
Second, suppose that $l^\nu(x^\nu,y^\nu) < -\rho$. Then, $f^\nu_{y^\nu}(\hat u^\nu,x^\nu) \leq -\rho\leq \alpha^\nu$ and we again have $(\hat u^\nu,x^\nu,\alpha^\nu)$ $\in$ $\epi f^\nu_{y^\nu} \cap \ball(\hat\rho)$. Thus, regardless of the case, there exist $(u,x,\alpha)\in \epi f_{y}$ such that $\max\{\|u - \hat u^\nu\|_2, \|x - x^\nu\|_2, |\alpha - \alpha^\nu| \} \leq \eta^\nu$,  where $\eta^\nu = \hatsetd_{\hat\rho}(\epi f_y, \epi f^\nu_{y^\nu}) + \epsilon$. Since $l(x,y)$ is defined in terms of minimization of $f_y(\,\cdot\,,x)$, one has $l(x,y) \leq f_{y}(u,x) \leq \alpha \leq \alpha^\nu + \eta^\nu$. Moreover, $\eta^\nu\in \reals$ implies that $(x,\alpha^\nu+\eta^\nu) \in \epi l(\,\cdot\,,y)$. The arbitrary choice of $(x^\nu,\alpha^\nu)$ yields
\[
\exs\big( \epi l^\nu(\,\cdot\,, y^\nu) \cap \ball(\rho); \epi l(\,\cdot\,,y)\big) \leq \eta^\nu.
\]
Since $\epsilon$ is arbitrary too, the left-hand side is bounded from above by $\hatsetd_{\hat\rho}(\epi f_y, \epi f^\nu_{y^\nu})$. Next, we repeat these arguments with the roles of $l^\nu(\,\cdot\,,y^\nu)$ and $l(\,\cdot\,,y)$ reversed, and invoke Theorem \ref{thm:errorRock}(b).\eop

\begin{figure}
\centering
\includegraphics[width=0.7\textwidth]{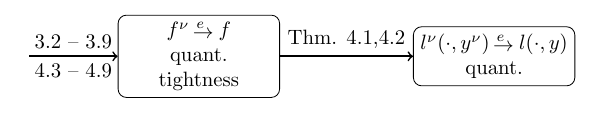}
\caption{\redrev{Usage of Theorems \ref{thm:epiLagrangians} and \ref{thm:errorLagrangians} to establish $l^\nu(\cdot,y^\nu)\eto l(\cdot,y)$ and its quantification.}}\label{overview4}
\end{figure}

\redrev{Figure \ref{overview4} shows how we can employ results from Section \ref{sec:tilted} and the present section to satisfy the assumptions of Theorems \ref{thm:epiLagrangians} and \ref{thm:errorLagrangians}.}  

\subsection{\redrev{Tightness in} Composite Optimization}\label{subsec:compositeLagr}

Theorem \ref{thm:epiLagrangians} brings forth the importance of tightness. We next examine sufficient conditions for tightness in the context of composite optimization as in Example \ref{e:composite}.

\begin{proposition}{\rm (composite optimization; tightness).}\label{prop:compositetight}
In the setting of Example \ref{e:composite}, suppose that the functions $\{h^\nu, \nu\in\nats\}$ are proper, lsc, and convex, and both $\{x^\nu, \nu\in\nats\}$ and $\{y^\nu, \nu\in\nats\}$ are convergent sequences. If $\{G^\nu(x^\nu), \nu\in\nats\}$ is bounded and there exists a compact set $K\subset\reals^m$ such that $\partial h^{\nu*}(y^\nu) \cap K \neq \emptyset$ for all $\nu\in\nats$, then $\{f_{y^\nu}^\nu(\,\cdot\,,x^\nu), \nu\in \nats\}$ is tight.
\end{proposition}
\state Proof. We can limit the focus to $x^\nu \in X^\nu$ because otherwise tightness holds trivially. Thus, without loss of generality, suppose that $x^\nu \in X^\nu$ for all $\nu\in \nats$. By \cite[Thm. 2.19]{primer} and \cite[Prop. 5.37]{primer}, \redrev{using the short-hand notation $w^\nu = G^\nu(x^\nu) + u^\nu$,} one has
\begin{equation}\label{eqn:unuequiv}
u^\nu \in \nargmin f_{y^\nu}^\nu(\,\cdot\,, x^\nu) ~\Longleftrightarrow~y^\nu \in \partial h^\nu(w^\nu) ~\Longleftrightarrow~ w^\nu \in \partial h^{\nu*}(y^\nu).
\end{equation}
By assumption, there exist points $\{v^\nu \in \partial h^{\nu*}(y^\nu) \cap K, \nu\in\nats\}$. Then, $\{u^\nu = v^\nu - G^\nu(x^\nu), \nu\in\nats\}$ is bounded because $v^\nu\in K$ and $\{G^\nu(x^\nu), \nu\in\nats\}$ is bounded. By \eqref{eqn:unuequiv}, $u^\nu \in \nargmin f_{y^\nu}^\nu(\,\cdot\,, x^\nu)$ and there is a compact set $B$ such that $\inf f_{y^\nu}^\nu(\,\cdot\,,x^\nu) = f_{y^\nu}^\nu(u^\nu,x^\nu)=\inf\{f_{y^\nu}^\nu(u,x^\nu)~|~ u\in B\}$ for all $\nu\in\nats$.\eop

The assumption about a compact set $K$ in Proposition \ref{prop:compositetight} holds broadly. For example, suppose that $h^\nu$ is of the form $h^\nu(z) = \sup_{y\in Y^\nu} \langle y,z\rangle - \tfrac{1}{2}\langle y, Q^\nu y\rangle$,
where $Y^\nu\subset\reals^m$ is nonempty, closed, and convex, and $Q^\nu$ is a symmetric positive semidefinite $m$-by-$m$ matrix; see, e.g., \cite[Sec. 5.I]{primer}. Then, for $y^\nu \in Y^\nu$, $\partial h^{\nu*}(y^\nu) = Q^\nu y^\nu + N_{Y^\nu}(y^\nu)$. Since $0 \in N_{Y^\nu}(y^\nu)$, the assumption of Proposition \ref{prop:compositetight} holds if $\{Q^\nu y^\nu, \nu\in\nats\}$ is bounded. This is indeed common. For example, the case $h^\nu = \iota_{(-\infty,0]^m}$ (i.e.,  inequality constraints) has  $Q^\nu = 0$ and $Y^\nu = [0,\infty)^m$. 


For fixed $y$, there is also a direct approach that does not rely on Theorem \ref{thm:epiLagrangians}.

\begin{proposition}{\rm (composite optimization; epi-convergence).}\label{prop:composite22}
In the setting of Example \ref{e:composite}, suppose that $X^\nu\sto X$ and $g_0^\nu(x^\nu)\to g_0(x)$ as well as $G^\nu(x^\nu) \to G(x)$ when $x^\nu\in X^\nu \to x\in X$. For $y\in \reals^m$, one has $l^\nu(\,\cdot\,, y) \eto l(\,\cdot\,,y)$ provided that $h^{\nu*}(y)\to h^*(y)$, which holds in particular when $y\not\in \bdry (\dom h^*)$, $h^\nu\eto h$, and $\{h, h^\nu, \nu\in\nats\}$ are proper, lsc, and convex functions.
\end{proposition}
\state Proof. Suppose that $h^{\nu*}(y)\to h^*(y)$. First, we consider the liminf-condition \eqref{eqn:liminf} and let $x^\nu\to x$. We observe that $\nliminf (\iota_{X^\nu}(x^\nu) - h^{\nu*}(y)) \geq \nliminf \iota_{X^\nu}(x^\nu)  + \nliminf (- h^{\nu*}(y))$ because $\iota_{X^\nu}(x^\nu)$ takes only the values $0$ and $\infty$. It then follows that
\[
\nliminf l^\nu(x^\nu,y) \geq \nliminf \iota_{X^\nu}(x^\nu) + \nliminf g_0^\nu(x^\nu) + \nliminf \hspace{-0.05cm}\big\langle G^\nu(x^\nu),y\big\rangle + \nliminf \hspace{-0.05cm}\big(- h^{\nu*}(y)\big).
\]
If $x\not\in X$, then the first term on the right-hand size is $\infty$ and $\nliminf l^\nu(x^\nu,y) \geq l(x,y)=\infty$. If $x\in X$, then we can assume without loss of generality that $x^\nu \in X^\nu$. This implies that $\nliminf g_0^\nu(x^\nu) = g_0(x)$ and $\nliminf \langle G^\nu(x^\nu),y\rangle = \langle G(x),y\rangle$. Thus, the liminf-condition \eqref{eqn:liminf} holds. Second, we consider the limsup-condition \eqref{eqn:limsup} and let $x\in \reals^n$. Without loss of generality, we assume that $x\in X$. Since $X^\nu\sto X$, there is $x^\nu\in X^\nu\to x$. Thus,
\begin{align*}
&\nlimsup l^\nu(x^\nu,y)\\
& \leq \nlimsup \iota_{X^\nu}(x^\nu) + \nlimsup g_0^\nu(x^\nu) + \nlimsup \big\langle G^\nu(x^\nu),y\big\rangle + \nlimsup \big(- h^{\nu*}(y)\big)\\
& = g_0(x) + \big\langle G(x),y\big\rangle - h^{*}(y) = l(x,y).
\end{align*}

For the claim about $y\not\in \bdry (\dom h^*)$, we observe that $h^{\nu*} \eto h^*$ by \cite[Wijsman's theorem 11.34]{VaAn} and we can invoke \cite[Thm. 7.17]{VaAn} to conclude that $h^{\nu*}(y) \to h^*(y)$.\eop

The condition $h^{\nu*}(y)\to h^*(y)$ in Proposition \ref{prop:composite22} holds trivially if $h^\nu = h$ because then $h^{\nu*} = h^*$.
A boundary point $y\in \bdry (\dom h^*)$ will have to be checked case by case. For proper, lsc, and convex $h$, the conjugate $h^*$ is real-valued if and only if $h$ is coercive \cite[Thm. 11.8(d)]{primer}. Roughly, a function is coercive if it grows faster than linear, with $\iota_{\{0\}^m}$ and $\half\theta \|\cdot\|_2^2$, $\theta>0$, being relevant examples. The function $h = \iota_{(-\infty,0]^m}$ is not coercive, with $h^* = \iota_{[0,\infty)^m}$ as its conjugate. Nevertheless, the approximation given by $h^\nu(z) = \nsum_{i=1}^m \theta^\nu \max\{0,z_i\}$ epi-converges to $h$ as $\theta^\nu\to \infty$. We have $h^{\nu*}(y) = \nsum_{i=1}^m \iota_{[0,\theta^\nu]}(y_i) = \iota_{[0,\theta^\nu]^m}(y)$. Thus, even on the boundary of $\dom h^*$, we have $h^{\nu*}(y) \to h^*(y)$.

The following result supports Theorem \ref{thm:errorLagrangians}; see Section \ref{sec:supp} for a proof.

\begin{proposition}{\rm (composite optimization; quantification).}\label{prop:compositequant}
In the setting of Example \ref{e:composite} for some fixed $\nu$, suppose that $h$ and $h^\nu$ are proper, lsc, and convex functions and there exists a compact set $K\subset\reals^m$ such that $\partial h^{\nu*}(y^\nu) \cap K \neq \emptyset$ and $\partial h^{*}(y) \cap K \neq \emptyset$. Let $\rho \in [0,\infty)$ and
\begin{align*}
\delta = &\max\Big\{\nsup_{x} \big\{ \big\|G(x)\big\|_2 ~\Big|~ \|x\|_2 \leq \rho, l(x,y) \leq \rho \big\},\\
&~~~~~~~~\nsup_{x} \big\{ \big\|G^\nu(x)\big\|_2 ~\Big|~ \|x\|_2 \leq \rho, l^\nu(x,y^\nu) \leq \rho \big\}\Big\}.
\end{align*}
Then, $\hat\rho \geq \max\{\rho, \delta + \sup_{v\in K} \|v\|_2\}$ satisfies the requirements in Theorem \ref{thm:errorLagrangians}.
\end{proposition}

\subsection{\redrev{Lagrangians for} Stochastic Optimization}

In this subsection, we examine the Lagrangians emerging from Propositions \ref{e:ambiguity} and \ref{e:splitting}.

\begin{example}{\rm (ambiguity).}\label{e:ambiguity2} In the setting of Proposition \ref{e:ambiguity}, $u_i \in [-p_i, 0]$ effectively for all $i = 1, \dots, m$ and the tightness assumption in Theorem \ref{thm:epiLagrangians} holds trivially. This establishes that the Lagrangians $l$ and $l^\nu$ produced by $f$ and $f^\nu$, respectively, convergence in the sense: $l^\nu(\,\cdot\,, y^\nu)\eto l(\,\cdot\,,y)$ whenever $y^\nu\to y$. Moreover, $\hat \rho$ can be set to $\max\{\rho, \|p\|_2\}$ in Theorem \ref{thm:errorLagrangians}. 

For $y\in \reals^m$ and $x\in \reals^n$, we obtain the following explicit formula using the definitions:  
\[
l(x,y) = g_0(x) + \nsum_{i=1}^m h_i(x), 
\]
where, for $\theta = 0$, $h_i(x) = p_i \min\{g_i(x), y_i\}$ and for $\theta \in (0,\infty)$,
\[
h_i(x) = \begin{cases}
    \tfrac{1}{2}\theta p_i^2 + y_i p_i & \mbox{ if } g_i(x) \in (y_i + \theta p_i, \infty]\\ 
    p_i g_i(x) + \tfrac{1}{2\theta}(y_i - g_i(x))^2 & \mbox{ if } g_i(x) \in [y_i, y_i + \theta p_i]\\ 
    p_i g_i(x) & \mbox{ if } g_i(x) \in (-\infty, y_i).
\end{cases}
\]
\end{example}

\begin{example}{\rm (splitting).}\label{e:splitting2}
In the setting of Proposition \ref{e:splitting}, suppose that $l$ and $l^\nu$ are the Lagrangians produced by $f$ and $f^\nu$, respectively. Then, we immediately obtain the explicit formula
\[
l(x,y) = - \nsum_{i=1}^m  p_i \star g_i^*(y_i) + \Big\langle x, \nsum_{i=1}^m y_i \Big\rangle, 
\]
with a parallel expression for $l^\nu$. Here, $\star$ specifies epi-multiplication; see \cite[p.24]{VaAn}. 

If $y^\nu\to y$, $p^\nu\to p$ with $p_i>0, i=1, \dots, m$, and $\{g_i^\nu, i=1, \dots, m, \nu\in\nats\}$ is equi-coercive (i.e., there is a coercive function $h$ such that $g_i^\nu(x) \geq h(x)$ for all $i$, $\nu$, $x$), then $l^\nu(\,\cdot\,, y^\nu) \eto l(\,\cdot\,,y)$. 
\end{example}
\state Detail. Epi-convergence follows from Theorem \ref{thm:epiLagrangians} after noting that tightness of $\{f_{y^\nu}^\nu(\,\cdot\,, x^\nu), \nu\in\nats\}$ holds if there is a bounded set containing a minimizer of $p_i^\nu g_i^\nu(x^\nu + u_i) - y_i^\nu u_i$ for each $i$ and $\nu$.\eop

\subsection{\redrev{Tightness from} Augmentation}

Formulas for Lagrangians under augmention are available in many cases; see \cite[Sec. 11.K]{VaAn}. Trivially, any Rockafellian $f:\reals^m\times\reals^n\to \Reals$ augmented with the indicator function $\iota_{\{0\}^m}$ produces the Rockafellian given by $\bar f(u,x) = f(u,x) + \iota_{\{0\}^m}(u)$, which in turn defines a Lagrangian of the form $\bar l(x,y) = f(0,x)$ for all $y$. We concentrate on two sufficient conditions for tightness as required in Theorem \ref{thm:epiLagrangians}. Augmentation is also discussed in \cite[Thm. 4.2]{AttouchAzeWets.86}, but for the convex case.

\begin{proposition}{\rm (tightness under proximal augmentation).}\label{e:tightAugment}
For $\theta \in (0,\infty)$, $\{\theta^\nu\in [\theta, \infty), \nu\in\nats\}$, and $\{f^\nu:\reals^m\times\reals^n\to \Reals, \nu\in\nats\}$, suppose that for each $\nu$, $f^\nu$ is a Rockafellian for $\phi^\nu:\reals^n\to \Reals$ and let 
\[
\bar f^\nu(u,x) = f^\nu(u,x) + \theta^\nu \|u\|^2_2. 
\]
Suppose that $\{x^\nu\in\reals^n, \nu\in\nats\}$ is convergent, $\{y^\nu\in\reals^m, \nu\in\nats\}$ is bounded, and there exist $\bar\nu\in \nats$, $\beta\in (0, \theta)$, and $\gamma, \tau\in \reals$ such that for each $\nu\geq \bar \nu$ with $\inf f^\nu(\,\cdot\,, x^\nu)<\infty$ one has 
\begin{align}
\forall u\in\reals^m: & ~~f^\nu(u,x^\nu) \geq \gamma - \beta \|u\|_2^2\label{eqn:conda2}\\ 
\exists u^\nu\in\reals^m  \mbox{ with } \|u^\nu\|_2 \leq \tau: & ~~f^\nu(u^\nu,x^\nu) \leq \tau\label{eqn:condb2}. 
\end{align}
Then, $\{\bar f^\nu_{y^\nu}(\,\cdot\,, x^\nu), \nu\in \nats\}$ is tight.
\end{proposition}
\state Proof. Since  $\{y^\nu, \nu\in\nats\}$ is bounded, there is $\rho < \infty$ such that $\sup_{\nu\in\nats} \|y^\nu\|_2\leq \rho$. Without loss of generality, we consider $\nu\geq \bar\nu$ with $\inf f^\nu(\,\cdot\,, x^\nu)<\infty$ because otherwise $\inf \bar f^\nu_{y^\nu}(\,\cdot\,,x^\nu) = \infty$ and the tigthness condition holds trivially. For such $\nu$, one has  
\begin{align*}
 \bar f^\nu_{y^\nu}(u,x^\nu) &= f^\nu(u,x^\nu) + \theta^\nu\|u\|_2^2 - \langle y^\nu, u\rangle \geq \gamma + (\theta^\nu - \beta) \|u\|_2^2 - \rho\|u\|_2~~~\forall u\in\reals^m\\
 \bar f^\nu_{y^\nu}(u^\nu,x^\nu) &= f^\nu(u^\nu,x^\nu) + \theta^\nu\|u^\nu\|_2^2 - \langle y^\nu, u^\nu\rangle \leq \tau + \theta^\nu \tau^2 + \rho \tau. 
\end{align*}
\redrev{Since $\theta^\nu - \beta>0$,} there exists $\bar\rho > \rho$, independent of $\nu$, such that $\gamma + (\theta^\nu - \beta) \|u\|_2^2 - \rho\|u\|_2 \geq \tau + \theta^\nu \tau^2 + \rho \tau$ when $\|u\|_2 > \bar\rho$. Consequently, $\bar f^\nu_{y^\nu}(u^\nu,x^\nu) \leq \bar f^\nu_{y^\nu}(u,x^\nu)$ when $\|u\|_2 > \bar\rho$ and $\nu\geq\bar\nu$. This implies that  
$\inf_{u\in \ball(\bar\rho)} \bar f^\nu_{y^\nu}(u,x^\nu)$ $\leq$ $\inf_{u\in \reals^m} \bar f^\nu_{y^\nu}(u,x^\nu)$ for all $\nu\geq \bar \nu$.\eop

\begin{proposition}{\rm (tightness under sharp augmentation).}\label{e:tightAugment2}
For $\{\theta^\nu\in (0, \infty), \nu\in\nats\}$ and $\{f^\nu:\reals^m\times\reals^n\to \Reals, \nu\in\nats\}$, suppose that for each $\nu$, $f^\nu$ is a Rockafellian for $\phi^\nu:\reals^n\to \Reals$ and let 
\[
\bar f^\nu(u,x) = f^\nu(u,x) + \theta^\nu \|u\|, 
\]
where the norm is arbitrary but equivalent to $\|\cdot\|_2$ in the sense that for some $\kappa \in (0,\infty)$ one has $\|u\|_2 \leq \kappa \|u\|$ for all $u\in\reals^m$. Suppose that $\{x^\nu\in\reals^n, \nu\in\nats\}$ is convergent, $\{y^\nu\in\reals^m, \nu\in\nats\}$ satisfies $\sup_{\nu\in\nats} \|y^\nu\|_2\leq \rho$ for some $\rho \in [0,\infty)$, and there exist $\bar\nu\in \nats$, $\beta\in (0, \infty)$, and $\gamma, \tau\in \reals$ such that for each $\nu\geq \bar \nu$ with $\inf f^\nu(\,\cdot\,, x^\nu)<\infty$ one has 
\begin{align}
\forall u\in\reals^m: & ~~f^\nu(u,x^\nu) \geq \gamma - \beta \|u\|\label{eqn:conda1}\\ 
\exists u^\nu\in\reals^m \mbox{ with } \|u^\nu\| \leq \tau: & ~~f^\nu(u^\nu,x^\nu) \leq \tau\label{eqn:condb1}. 
\end{align}
If $\theta^\nu \geq \beta + \kappa\rho + 1$ for all $\nu\geq \bar\nu$, then $\{\bar f^\nu_{y^\nu}(\,\cdot\,, x^\nu), \nu\in \nats\}$ is tight.
\end{proposition}
\state Proof. The conclusion follows by arguments similar to those for Proposition \ref{e:tightAugment}.\eop

The conditions \eqref{eqn:conda2} and \eqref{eqn:conda1} are mild as they simply impose lower bounding concave functions on the Rockafellians $f^\nu$ uniformly for all $\nu$. Conditions \eqref{eqn:condb2} and \eqref{eqn:condb1} need to be verified on a case-by-case basis, but hold, for example, in the composite setting of Proposition \ref{prop:composite}.

The path to analyzing many Lagrangians constructed by augmentation is now clear. For example, consider the Rockafellians $f$ and $f^\nu$ in Example \ref{e:composite} for composite optimization. Proposition \ref{prop:composite} furnishes sufficient conditions for $f^\nu\eto f$. 
Suppose that $f$ is augmented with $\iota_{\{0\}^m}$ to produce $\bar f$, which then is strictly exact, and $f^\nu$ is augmented by some norm to produce $\bar f^\nu$ as in Proposition \ref{prop:compositeAug}. That proposition gives sufficient conditions for $\bar f^\nu \eto \bar f$. This in turn feeds into Theorem \ref{thm:epiLagrangians}, which we can use to establish that the corresponding Lagrangians satisfy $\bar l^\nu(\,\cdot\,, y^\nu)\eto \bar l(\,\cdot\,, y)$. The theorem requires a tightness condition to hold, and this is where Propositions \ref{e:tightAugment} and \ref{e:tightAugment2} enter; \redrev{see Figure \ref{overview4}}.

\section{Convergence of Dual Functions}

Next, we focus on hypo-convergence of $\psi^\nu$ to $\psi$ and its quantification as represented by the top arrow in Figure \ref{overview}. There are three possible perspectives when examining a dual problem of maximizing $\psi$ and its approximation involving $\psi^\nu$: 
(i) One could develop sufficient conditions for $\psi^\nu \hto \psi$ expressed in terms of the underlying Rockafellians. (ii) One could focus on the corresponding Lagrangians and the max-inf problems $\nnmax_{y\in \reals^m} \ninf_{x\in \reals^n} l(x,y)$ and $\nnmax_{y\in \reals^m} \ninf_{x\in \reals^n} l^\nu(x,y)$. (iii) One could focus on the bifunctions $k,k^\nu:\reals^m\times \reals^{m+n}\to \Reals$ given by 
\begin{equation}\label{eqn:kbifcns}
k\big(y,(u,x)\big) =  f(u,x) - \langle y,u\rangle~~~\mbox{ and } ~~~k^\nu\big(y,(u,x)\big) =  f^\nu(u,x) - \langle y,u\rangle, 
\end{equation} 
and the resulting max-inf problems 
\[
\nnmax_{y\in \reals^m} \ninf_{(u,x)\in \reals^{m+n}} k\big(y,(u,x)\big) ~~\mbox{ and } ~~ \nnmax_{y\in \reals^m} \ninf_{(u,x)\in \reals^{m+n}} k^\nu\big(y,(u,x)\big). 
\]

We adopt perspective (i) in the following because it directly addresses the dual problems without passing through some equivalent (in some sense) max-inf problems. Perspectives (ii) and (iii) land us in the territory of lopsided convergence as pioneered in \cite{AttouchWets.83b} and refined in \cite{RoysetWets.19a}. That general approach to analyzing max-inf (or min-sup) problems furnishes sufficient conditions for the convergence of max-inf points and max-inf values. However, the conditions do not explicitly leverage the special structure of the bifunctions $l,l^\nu$ and $k,k^\nu$, in particular concavity in the $y$-arguments. Still, we compare below with (iii) in detail. Perspective (ii) would result in conditions on the Lagrangians, which is counter to our focus on the more fundamental Rockafellians.  

The next theorem provides the key building block for examining dual functions in applications. 

\begin{theorem}{\rm (hypo-convergence of dual functions).}\label{thm:hypoofdual}
For Rockafellians $f,f^\nu:\reals^m\times\reals^n\to \Reals$, suppose that $f^\nu \eto f$. Then, the following hold:
\begin{enumerate}[(a)]

\item If for each $y\in \dom (-\psi)$ there exists $y^\nu \to y$ such that $\{f^\nu_{y^\nu}, \nu\in\nats\}$ is tight, then $\psi^\nu \hto \psi$.
In fact, for all such $y^\nu \to y$, one has $\psi^\nu(y^\nu) \to \psi(y)$.

\item  If $\dom(-\psi)$ has a nonempty interior, $Y$ is a dense subset of $\dom (-\psi)$, and $\{f^\nu_{y}, \nu\in\nats\}$ is tight for each $y\in Y$, then $\psi^\nu \hto \psi$.

\end{enumerate}

\end{theorem}
\state Proof. For (a), it suffices to show that $-\psi^\nu \eto -\hspace{-0.035cm}\psi$. First, we consider \eqref{eqn:liminf} and let $y^\nu\to y$. By Theorem \ref{thm:errorRock}(a), $f^\nu_{y^\nu} \eto f_y$. In view of Proposition \ref{prop:tigthepi}(b), this implies that $\nlimsup ( \ninf f^\nu_{y^\nu} ) \leq \ninf  f_y$. Thus, $\nliminf -\psi^\nu(y^\nu) \geq -\psi(y)$. Second, we consider \eqref{eqn:limsup}. Let $y\in \reals^m$. If $y \not\in \dom (-\psi)$, then $\nlimsup -\psi^\nu(y) \leq -\psi(y)$ holds trivially because the right-hand side is infinity. Thus, we concentrate on $y\in \dom (-\psi)$. By assumption, there exists $y^\nu \to y$ such that $\{f^\nu_{y^\nu}, \nu\in\nats\}$ is tight. By Theorem \ref{thm:errorRock}(a), $f^\nu_{y^\nu} \eto f_y$. Proposition \ref{prop:tigthepi}(c) then implies that \begin{equation}\label{eqn:infconv}
\inf  f^\nu_{y^\nu} \to \inf  f_y.
\end{equation}
Thus, $\nlimsup -\psi^\nu(y^\nu) \leq -\psi(y)$ and the first part of (a) follows. Since \eqref{eqn:infconv} simply states that $\psi^\nu(y^\nu) \to \psi(y)$, the remaining claim in (a) holds immediately.

For (b), we start by showing that $-\psi$ is lsc. We consider two cases. (i) Suppose that $f(u,x) = -\infty$ for some $(u,x)\in \reals^m\times\reals^n$. Then, $-\psi(y) = \infty$ for all $y\in\reals^m$. This means that $\dom(-\psi) = \emptyset$, which is ruled out by assumption. (ii) Suppose that $f(u,x)>-\infty$ for all $(u,x)\in \reals^m\times\reals^n$. Since $-\psi$ can be written as $-\psi(y) = \nsup\{ \langle y,u\rangle - f(u,x) | (u,x) \in \dom f \}$ for all $y\in \reals^m$, we invoke \cite[Prop. 6.27]{primer} to conclude that $-\psi$ is lsc. That proposition requires $y\mapsto \langle y,u\rangle - f(u,x)$ to be real-valued and lsc for each $(u,x) \in \dom f$, which indeed hold.

Since $-\psi^\nu$ and $-\psi$ are convex functions with $\dom(-\psi)$ having a nonempty interior and $-\psi$ being lsc, it suffices by \cite[Thm. 7.17]{primer} to show that $\psi^\nu(y) \to \psi(y)$ for all $y\in Y \cup (\reals^m \setminus \dom (-\psi))$. Let $y\in Y \cup (\reals^m \setminus \dom (-\psi))$.

We invoke Theorem \ref{thm:errorRock}(a) to establish that $f^\nu_{y} \eto f_y$ and consider two cases. (i) Suppose that $y\in Y$. Then, $\{f^\nu_{y}, \nu\in\nats\}$ is tight and Proposition \ref{prop:tigthepi}(c) implies that $\inf  f^\nu_{y} \to \inf  f_y$, which is identical to having $\psi^\nu(y) \to \psi(y)$. (ii) Suppose that $y\not\in \dom (-\psi)$. Then, $\psi(y) = -\infty$, which is equivalent to having $\inf f_y = -\infty$. By Proposition \ref{prop:tigthepi}(b), $\nlimsup (\inf f_y^\nu) \leq \inf f_y = -\infty$ and $\inf f_y^\nu \to \inf f_y$.\eop

Comparing the assumption in Theorem \ref{thm:hypoofdual}(a) with that of Theorem \ref{thm:hypoofdual}(b), we see that the former requires a check for each $y\in \dom(-\psi)$ while the latter only considers a dense subset of such points. But, the former permits $y^\nu\neq y$ while the latter needs $y^\nu = y$ in the tightness verification. 

The connection between Theorem \ref{thm:hypoofdual} and perspective (iii), with bifunctions $k^\nu,k$ defined in \eqref{eqn:kbifcns}, is as follows. By \cite[Propostion 9.28]{primer}, if $k^\nu$ lop-converges to $k$ ancillary-tightly, then $\psi^\nu \hto \psi$. A close examination of the respective definitions (cf. Def. 9.14, 9.19, and 9.21 in \cite{primer}) reveals that $f^\nu\eto f$ if and only if $k^\nu$ lop-converges to $k$. However, the tightness assumptions in Theorem \ref{thm:hypoofdual}(a,b) are slightly weaker than ancillary-tightness, which requires examining all $y\in \reals^m$; see \cite[Def. 9.21]{primer}. The direct proof of Theorem \ref{thm:hypoofdual} has the merit that it completely bypasses the intricacies of lopsided convergence. Theorem \ref{thm:hypoofdual}(a) also furnishes additional insight about continuous convergence.

The property $\psi^\nu \hto \psi$ together with tightness of $\{-\psi^\nu, \nu\in\nats\}$ ensure via Proposition \ref{prop:tigthepi}(c) that $\sup \psi^\nu \to \sup\psi$. If $\sup \psi = \inf f(0, \,\cdot\,)$ also holds, then we conclude that the maximum values of the approximating dual functions tend to the minimum value of the actual problem.

It is well-known that maximizers of the dual function $\psi$ are subgradients of the min-value function $u\mapsto \inf f(u, \,\cdot\,)$ under certain assumptions that typically include convexity of the underlying Rockafellian $f$; see, e.g., \cite[Thm. 5.44]{primer}. The property $\psi^\nu \hto \psi$, therefore, helps us to confirm that maximizers of $\psi^\nu$ can only converge to such subgradients.

Even in the convex case, tightness is hard to avoid. If $f,f^\nu:\reals^m\times\reals^n\to \Reals$ are proper, lsc, and convex Rockafellians, then Wijsman's theorem (see, e.g., \cite[Thm. 11.34]{VaAn}) ensures that $f^\nu \eto f$ if and only if $f^{\nu*}\eto f^*$. Since $-\psi(y) = f^*(y,0)$ and $-\psi^\nu(y^\nu) = f^{\nu*}(y^\nu,0)$, we have as a consequence of $f^\nu \eto f$ that $\nliminf -\psi^\nu(y^\nu) \geq -\psi(y)$ whenever $y^\nu\to y$. This provides ``one-half'' of the requirement for $-\psi^\nu\eto -\psi$; see \eqref{eqn:liminf}. For the other ``half'' \eqref{eqn:limsup}, we seek for each $y$ a sequence $y^\nu\to y$ such that $\nlimsup -\psi^{\nu}(y^\nu) \leq -\psi(y)$ or equivalently $\nlimsup f^{\nu*}(y^\nu,0) \leq f^*(y,0)$, which is identical to $\nliminf (\inf f^\nu_{y^\nu} ) \geq \inf f_y$. We are again compelled to bring in tightness as in Proposition \ref{prop:tigthepi}(a).

The truncated Hausdorff distance between $\epi f$ and $\epi f^\nu$ also leads to a bound for the dual problems. We note that $\hatsetd_\rho(\hypo \psi, \hypo \psi^\nu)$ directly quantifies the distance between maximum values and near-maximizers as we see from Theorem \ref{tapproxoptimalvalue} after a re-orientation from minimization to maximization. 

The bifunction perspective (iii) discussed above can also lead to estimates of $\hatsetd_\rho(\hypo \psi^\nu, \hypo \psi)$ as seen in \cite[Sec. 5]{RoysetWets.16b}. However, these estimates rely on the pointwise difference between $k^\nu$ and $k$, Lipschitz properties, and also fail to utilize the specific structure of these bifunctions. Thus, the following bound is tighter and more versatile.

\begin{theorem}{\rm (error bounds for dual functions).}\label{thm:quantDual}
For Rockafellians $f,f^\nu:\reals^m \times \reals^n\to \Reals$, suppose that $\epi f$, $\epi f^\nu$, $\dom(-\psi)$, and $\dom(-\psi^\nu)$ are nonempty, and $\rho \in [0,\infty)$. Let $\|\cdot\|_2$, $(u,x)\mapsto \max\{\|u\|_2, \|x\|_2\}$, $(y,\alpha)\mapsto \max\{\|y\|_2, |\alpha|\}$, and $(u,x,\alpha)\mapsto \max\{\|u\|_2, \|x\|_2, |\alpha|\}$ be the norms on $\reals^m$, $\reals^m\times \reals^n$, $\reals^m\times \reals$, and $\reals^m\times\reals^n\times \reals$, respectively. Then, the following hold:
\begin{enumerate}[(a)]

\item If $\hat\rho > \rho + \hatsetd_\rho(\dom(-\psi), \dom(-\psi^\nu))$ and $\rho'\in [1,\infty)$ satisfies
\begin{align*}
  &\forall y\in \dom(-\psi) \cap \ball(\hat\rho): &&\nargmin f_y \cap \ball(\rho')\neq \emptyset, ~~~\psi(y) \in [-\rho',\rho']\\
  &\forall y^\nu\in \dom(-\psi^\nu) \cap \ball(\hat\rho): &&\nargmin f^\nu_{y^\nu} \cap \ball(\rho')\neq \emptyset, ~~~\psi^\nu(y^\nu) \in [-\rho',\rho'],
\end{align*}
then, provided that $\bar\rho \geq \rho'(1 + \hat\rho)$, one has
\[
\hatsetd_\rho(\hypo \psi^\nu, \hypo \psi) \leq (1+\hat\rho)\hatsetd_{\bar\rho}(\epi f, \epi f^\nu) + \rho'\hatsetd_\rho\big(\dom(-\psi), \dom(-\psi^\nu) \big).
\]

\item If $Y(\rho) = \dom(-\psi) \cap \dom(-\psi^\nu) \cap \ball(\rho)$ is nonempty and $\rho'\in [0,\infty)$ satisfies for all $y\in Y(\rho)$
\[
\nargmin f_y \cap \ball(\rho')\neq \emptyset, ~~~\nargmin f^\nu_{y} \cap \ball(\rho')\neq \emptyset, ~~~\psi(y), \psi^\nu(y) \in [-\rho',\rho'],
\]
then, provided that $\bar\rho \geq \rho'(1 + \rho)$, one has 
\[
\sup_{y\in Y(\rho)} \big|\psi(y) - \psi^\nu(y)\big| \leq (1+\rho)\hatsetd_{\bar\rho}(\epi f, \epi f^\nu).
\]

\end{enumerate}

\end{theorem}
\state Proof. For (a), let $(y^\nu,\alpha^\nu) \in \hypo \psi^\nu \cap \ball(\rho)$. Then, $y^\nu \in \dom(-\psi^\nu)$. We consider two cases.
First, suppose that $y^\nu\in \dom(-\psi)$. Set $y = y^\nu$. Since $\hat\rho>\rho$, the assumptions establish that
\begin{equation}\label{eqn:dualerror1}
\nargmin f_y \cap \ball(\rho') \neq \emptyset, ~~~\nargmin f^\nu_{y^\nu} \cap \ball(\rho') \neq \emptyset, ~~~ \inf f_y, \inf f_{y^\nu}^\nu \in [-\rho',\rho'].
\end{equation}
We can therefore apply Proposition \ref{tapproxoptimalvalue}(a) and establish
\[
\big|\psi(y) - \psi^\nu(y^\nu)\big| = |\inf f_y - \inf f^\nu_{y^\nu}| \leq \hatsetd_{\rho'}(\epi f_y, \epi f^\nu_{y^\nu}) \leq (1+\rho)\hatsetd_{\breve\rho}(\epi f, \epi f^\nu),
\]
with $\breve\rho \geq \rho'(1+\rho)$, where the last inequality follows from Theorem \ref{thm:errorRock}(b). Thus, $(y,\min\{\alpha^\nu,\psi(y)\}) \in \hypo \psi$ and $|\alpha^\nu - \min\{\alpha^\nu,\psi(y)\}| \leq (1+\rho)\hatsetd_{\breve\rho}(\epi f, \epi f^\nu)$. This means that the distance between $(y^\nu,\alpha^\nu)$ and $(y,\min\{\alpha^\nu,\psi(y)\})$ in the chosen norm is at most $(1+\rho)\hatsetd_{\breve\rho}(\epi f, \epi f^\nu)$.

Second, suppose that $y^\nu\hspace{-0.04cm}\not\in \dom(-\psi)$. Let $\epsilon \in (0, \hat\rho - \rho - \hatsetd_\rho(\dom(-\psi), \dom(-\psi^\nu) )]$. Since $y^\nu \in \dom(-\psi^\nu) \cap \ball(\rho)$, there is $y\in \dom(-\psi)$ such that 
\[
\|y- y^\nu\|_2 \leq \hatsetd_\rho\big(\dom(-\psi), \dom(-\psi^\nu) \big) + \epsilon.
\]
Thus,
\[
\|y\|_2 \leq \|y^\nu\|_2 + \|y- y^\nu\|_2 \leq \rho + \hatsetd_\rho\big(\dom(-\psi), \dom(-\psi^\nu) \big) + \epsilon \leq \hat\rho.
\]
By assumption, \eqref{eqn:dualerror1} again holds and we can bring in Proposition \ref{tapproxoptimalvalue}(a) and Theorem \ref{thm:errorRock}(b) to establish
\begin{align*}
\big|\psi(y) - \psi^\nu(y^\nu)\big| = |\inf f_y - \inf f^\nu_{y^\nu}| & \leq \hatsetd_{\rho'}(\epi f_y, \epi f^\nu_{y^\nu})\\
& \leq (1+\hat\rho)\hatsetd_{\bar\rho}(\epi f, \epi f^\nu) + \rho'\|y - y^\nu\|_2,
\end{align*}
with $\bar\rho \geq \rho'(1+\hat\rho)$. Thus, $(y,\min\{\alpha^\nu,\psi(y)\}) \in \hypo \psi$ and, with 
\[
\eta_\epsilon = (1+\hat\rho)\hatsetd_{\bar\rho}(\epi f, \epi f^\nu) + \rho'\hatsetd_\rho\big(\dom(-\psi), \dom(-\psi^\nu) \big) + \rho'\epsilon,
\]
one has 
\[
\big|\alpha^\nu - \min\{\alpha^\nu,\psi(y)\}\big| \leq (1+\hat\rho)\hatsetd_{\bar\rho}(\epi f, \epi f^\nu) + \rho'\|y - y^\nu\|_2 \leq \eta_\epsilon.
\]
The norm of $(y^\nu,\alpha^\nu)-(y,\min\{\alpha^\nu,\psi(y)\})$ is then at most
\[
\max\Big\{\hatsetd_\rho\big(\dom(-\psi), \dom(-\psi^\nu) \big) +\epsilon, ~\eta_\epsilon\Big\}.
\] 

Since $\rho'\geq 1$, $\hat\rho \geq \rho$, $\bar\rho\geq \breve\rho$, and $\hatsetd_{\bar\rho}(\epi f, \epi f^\nu) \geq \hatsetd_{\breve\rho}(\epi f, \epi f^\nu)$, we find that, in both cases, the distance between $(y^\nu,\alpha^\nu)$ and $(y,\min\{\alpha^\nu,\psi(y)\})$ in the chosen norm is at most $\eta_\epsilon$.
Moreover, $(y^\nu,\alpha^\nu)$ is arbitrary and this yields $\exs(\hypo \psi^\nu \cap \ball(\rho); \hypo \psi) \leq \eta_\epsilon$. We repeat this argument with the roles of $\psi$ and $\psi^\nu$ reversed and obtain $\hatsetd_\rho(\hypo \psi^\nu, \hypo \psi) \leq \eta_\epsilon$. Since $\epsilon$ is arbitrary, the claim holds.

For (b), let $y\in Y(\rho)$. The theorem's assumptions establish that $\nargmin f_y \cap \ball(\rho')$ and $\nargmin f^\nu_{y} \cap \ball(\rho')$ are nonempty and that $\inf f_y, \inf f_{y}^\nu \in [-\rho',\rho']$. By Proposition \ref{tapproxoptimalvalue}(a), this ensures that
\[
\big|\psi(y) - \psi^\nu(y)\big| = |\inf f_y - \inf f^\nu_{y}| \leq \hatsetd_{\rho'}(\epi f_y, \epi f^\nu_{y}) \leq (1+\rho)\hatsetd_{\bar\rho}(\epi f, \epi f^\nu),
\]
with $\bar\rho \geq \rho'(1+\rho)$, where the last inequality follows from Theorem \ref{thm:errorRock}(b).\eop

\begin{figure}
\centering
\includegraphics[width=0.7\textwidth]{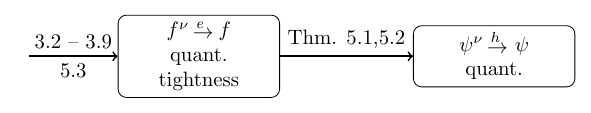}
\caption{\redrev{Usage of Theorems \ref{thm:hypoofdual} and \ref{thm:quantDual} to establish $\psi^\nu\hto\psi$ and its quantification.}}\label{overview5}
\end{figure}

\redrev{Figure \ref{overview5} gives an overview of how Theorems \ref{thm:hypoofdual} and \ref{thm:quantDual} can employ results from Section \ref{sec:tilted} and the present section to satisfy their assumptions.}  

\subsection{\redrev{Tightness in} Composite Optimization}\label{subsec:compositeDual}

In composite optimization, the tightness assumption of Theorem \ref{thm:hypoofdual} is addressed by the following facts. 

\begin{proposition}{\rm (composite optimization; tightness).}\label{prop:compositetight2}
In the setting of Example \ref{e:composite}, suppose that $\{h,h^\nu, \nu\in\nats\}$ are proper, lsc, and convex functions with $h^\nu\eto h$, $\cup_{\nu\in\nats} X^\nu$ is a bounded set, and $\sup\{\|G^\nu(x^\nu)\|_2\,|\, x^\nu\in X^\nu, \nu\in\nats\}$ is finite. If $y \in \reals^m$, then the following hold:
\begin{enumerate}[(a)]

\item If $\partial h^*(y)\neq \emptyset$, then there is $y^\nu\to y$ such that $\{f_{y^\nu}^\nu, \nu\in \nats\}$ is tight.

\item If there exist a compact set $K\subset\reals^m$ and a positive integer $\bar\nu$ such that $\partial h^{\nu*}(y) \cap K \neq \emptyset$ for all $\nu\geq \bar\nu$, then $\{f_{y}^\nu, \nu\in \nats\}$ is tight.

\end{enumerate}
\end{proposition}
\state Proof. For (a), suppose that $y$ is as stipulated in the proposition. Since $h^{\nu*}\eto h^*$ by \cite[Wijsman's theorem 11.34]{VaAn}, it follows from \cite[Attouch's theorem 12.35]{VaAn} that $\{(y,v) \, | \, v \in \partial h^{\nu*}(y)\}\sto \{(y,v) \, | \, v \in \partial h^{*}(y)\}$. By assumption, there exists $v\in \partial h^*(y)$ and thus also $v^\nu\to v$ and $y^\nu\to y$ such that $v^\nu \in \partial h^{\nu*}(y^\nu)$. Moreover, there is a bounded set $D\subset\reals^m$ such that $G^\nu(x^\nu) \in D$ for all $x^\nu \in X^\nu$ and $\nu\in\nats$. Let $\rho = \sup_{\nu\in\nats} \|v^\nu\|_2$, which is finite because $\{v^\nu, \nu\in\nats\}$ is a convergent sequence. We construct the compact set $B = \{v'-w\in\reals^m~|~\|v'\|_2\leq \rho, ~w\in D\}$.
Let $\{x^\nu\in \reals^n, \nu\in\nats\}$ be an arbitrary sequence. We consider two cases.

First, suppose that $x^\nu\in X^\nu$ for all $\nu\in\nats$. Then, we construct $u^\nu = v^\nu - G^\nu(x^\nu)$. Consequently, $u^\nu \in \nargmin f_{y^\nu}^\nu(\,\cdot\,, x^\nu)$ by \eqref{eqn:unuequiv}. Moreover, $\{u^\nu, \nu\in\nats\}\subset B$. This means that
\begin{equation}\label{eqn:tightnesscheckinf}
\inf\big\{f_{y^\nu}^\nu(u,x^\nu)~\big|~ u\in B\big\} \leq \inf f_{y^\nu}^\nu(\,\cdot\,,x^\nu) ~~~~\forall \nu\in\nats.
\end{equation}

Second, suppose that $x^\nu\not\in X^\nu$ for one or more $\nu\in\nats$. For such $\nu$, $f_{y^\nu}^\nu(u,x^\nu) = \infty$ regardless of $u\in\reals^m$. Thus, \eqref{eqn:tightnesscheckinf} still holds because each side of that inequality is $\infty$ for such $\nu$.

Let $\bar X = \cl \cup_{\nu\in \nats} X^\nu$ and fix $\nu$ in \eqref{eqn:tightnesscheckinf}. Since $x^\nu$ is arbitrary, we obtain from that inequality the relation $\inf\{f_{y^\nu}^\nu(u,x)~|~ u\in B, x\in \bar X\} \leq \inf f_{y^\nu}^\nu$. Since $\nu$ is arbitrary, this proves the claim.

For (b), suppose that $y$, $K$, and $\bar\nu$ are as stipulated in the proposition. By assumption, there is $\{v^\nu \in \partial h^{\nu*}(y) \cap K, \nu\geq \bar \nu\}$. We can then mimic the proof of (a) and reach the conclusion.\eop

Proposition \ref{prop:compositetight2}(a) helps to confirm the assumption of Theorem \ref{thm:hypoofdual}(a) because when $X$ is compact, $g_0$ is lsc, and $G$ is continuous, then $\dom(-\psi) = \dom h^*$.
The condition $\partial h^{*}(y)\neq \emptyset$ holds always when $y \in \nt(\dom h^*)$, but also on the whole domain of $h^*$ in common cases. It can be checked using any of the three equivalent conditions: $z \in \nargmin h - \langle y, \cdot\,\rangle$, $y \in \partial h(z)$, or $z \in \partial h^{*}(y)$.

Proposition \ref{prop:compositetight2}(b) supports Theorem \ref{thm:hypoofdual}(b). For example, if $h(z) = \iota_{(-\infty,0]^m}(z)$  and $h^\nu(z)$ $=$ $\theta^\nu \nsum_{i=1}^m \max\{0,z_i\}$, with $\theta^\nu \in (0,\infty)$, then we find that $\partial h^{\nu*}(y) = N_{[0,\theta^\nu]^m}(y)$. Under the assumption that $\dom(-\psi) = \dom h^*$ and $\theta^\nu\to \infty$, we see that $y\in \dom(-\psi)$ implies that $\partial h^{\nu*}(y)$ is nonempty and contains the zero vector for sufficiently large $\nu$. Thus, one can take $K = \{0\}^m$ in Proposition \ref{prop:compositetight2}(b).

\section{Additional Material and Proofs}\label{sec:supp}

The consequences of epi-convergence are well known. The following facts can be proven using arguments nearly identical to those underpinning \cite[Thm. 5.5]{primer}.

\begin{proposition}{\rm (consequences of epi-convergence).}\label{prop:tigthepi}
For $\{g, g^\nu:\reals^n\to \Reals, \nu\in\nats\}$, one has:
\begin{enumerate}[(a)]

\item If  $\nliminf g^\nu(x^\nu) \geq g(x)$ whenever $x^\nu \to x$ and $\{g^\nu, \nu\in\nats\}$ is tight, then  
\[
\nliminf (\inf g^\nu) \geq \inf g.
\]

\item If for all $x\in \reals^n$ there is $x^\nu\to x$ such that $\nlimsup g^\nu(x^\nu) \leq g(x)$, then 
\[
\nlimsup (\inf g^\nu) \leq \inf g.
\]

\item If $g^\nu\eto g$ and $\{g^\nu, \nu\in\nats\}$ is tight, then 
\[
\inf g^\nu \to \inf g.
\]

\item If $g^\nu\eto g$, $\inf g < \infty$, $x^\nu \in \epsilon^\nu\mbox{-}\argmin g^\nu$, $\epsilon^\nu\downto 0$, and $x^\nu\Nto x$ for some subsequence $N\subset \nats$, then 
    \[
    x\in \nargmin g ~\mbox{ and }~ g^\nu(x^\nu) \Nto \inf g.
    \]
    
\end{enumerate}
\end{proposition}

\state Proof of Theorem \ref{thm:errorRock}. For (a), since $(u,x)\mapsto \langle y^\nu, u\rangle$ both epi-converges and hypo-converges to $(u,x)\mapsto \langle y, u\rangle$ and the latter function is real-valued, we can follow the proof of \cite[Prop. 4.19]{primer}.

For (b), let $\epsilon \in (0,\infty)$ and $(u^\nu,x^\nu,\alpha^\nu) \in \epi f^\nu_{y^\nu} \cap \ball(\rho)$. Thus,
$f^\nu(u^\nu,x^\nu) \leq \alpha^\nu + \langle y^\nu, u^\nu\rangle$ and
\[
\big|\alpha^\nu+ \langle y^\nu, u^\nu\rangle\big| \leq |\alpha^\nu| + \|y^\nu\|_2\|u^\nu\|_2 \leq \rho\big(1 + \|y^\nu\|_2\big)\leq \bar\rho.
\]
Since $(u^\nu,x^\nu,\alpha^\nu+ \langle y^\nu, u^\nu\rangle) \in \epi f^\nu \cap \ball(\bar\rho)$, there exist $(u,x,\alpha)\in \epi f$ such that
\[
\max\big\{\|u-u^\nu\|_2, \|x-x^\nu\|_2, \big|\alpha-\alpha^\nu-\langle y^\nu,u^\nu\rangle\big|\big\}  \leq \hatsetd_{\bar\rho}(\epi f, \epi f^\nu) + \epsilon.
\]
With $\delta = \hatsetd_{\bar\rho}(\epi f, \epi f^\nu) + \epsilon + \rho\|y-y^\nu\|_2 + \|y\|_2(\hatsetd_{\bar\rho}(\epi f, \epi f^\nu) + \epsilon)$, the point $(u,x)$ satisfies
\[
f_y(u,x) \leq \alpha - \langle y, u\rangle \leq \alpha^\nu + \langle y^\nu,u^\nu\rangle + \hatsetd_{\bar\rho}(\epi f, \epi f^\nu) + \epsilon - \langle y, u\rangle \leq \alpha^\nu + \delta.
\]
Since $\delta\in \reals$, this means that $(u,x,\alpha^\nu+\delta) \in \epi f_y$ and its distance to $(u^\nu,x^\nu,\alpha^\nu)$ in the adopted norm is at most $\delta$. Since $(u^\nu,x^\nu,\alpha^\nu) \in \epi f^\nu_{y^\nu} \cap \ball(\rho)$ is arbitrary, we find that
$\exs( \epi f^\nu_{y^\nu} \cap \ball(\rho);  \epi f_y  ) \leq \delta$. Since $\epsilon$ is arbitrary, it also follows that
\[
\exs\big( \epi f^\nu_{y^\nu} \cap \ball(\rho);  \epi f_y  \big) \leq \Big( 1 + \max\big\{\|y\|_2, \|y^\nu\|_2\big\} \Big)\hatsetd_{\bar\rho}(\epi f, \epi f^\nu) + \rho \|y-y^\nu\|_2.
\]
After repeating the arguments with the roles of $f_y$ and $f^\nu_{y^\nu}$ reversed, we reach the conclusion.\eop

\state Proof of Proposition \ref{prop:composite}. We first consider the liminf-condition \eqref{eqn:liminf}. Let $u^\nu\to u$ and $x^\nu\to x$. Since $X^\nu\sto X$, $\nliminf \iota_{X^\nu}(x^\nu) \geq \iota_X(x)$. If $x\not\in X$, then this implies that $\nliminf f^\nu(u^\nu,x^\nu) \geq f(u,x)$ because $h$ is proper. If $x\in X$, then we proceed without loss generality under the assumption that $x^\nu \in X^\nu$ because $f^\nu(u^\nu,x^\nu) = \infty$ otherwise. Thus, $g_0^\nu(x^\nu)\to g_0(x)$ and $G^\nu(x^\nu)\to G(x)$. Since $h^\nu\eto h$, one has $\nliminf h^\nu(G^\nu(x^\nu) + u^\nu)\geq h(G(x) + u)$. We conclude that $\nliminf f^\nu(u^\nu,x^\nu) \geq f(u,x)$ holds in this case too. Second, we consider the limsup-condition \eqref{eqn:limsup}. Let $(u,x)\in \reals^m\times \reals^n$. Without loss of generality, we assume that $x\in X$. Since $X^\nu\sto X$, there is $x^\nu\in X^\nu\to x$. Also, since $h^\nu\eto h$, there is $z^\nu$ $\to$ $G(x) + u$ such that $\nlimsup h^\nu(z^\nu) \leq h(G(x) + u)$. Set $u^\nu = z^\nu-G^\nu(x^\nu)$. Then, $u^\nu\to u$. Moreover,
\begin{align*}
\nlimsup f^\nu(u^\nu,x^\nu) & \leq \nlimsup g_0^\nu(x^\nu) + \nlimsup h^\nu\big(G^\nu(x^\nu) + u^\nu\big)\\
& \leq g_0(x) + h\big(G(x) + u\big) = f(u,x)
\end{align*}
and the conclusion follows.\eop

\state Proof of Theorem \ref{thm:augmentation}. Trivially, $\bar f$ and $\bar f^\nu$ are Rockafellians for $\phi$ and $\phi^\nu$. The liminf- and limsup-conditions in \eqref{eqn:liminf} and \eqref{eqn:limsup} immediately lead to (a). For (b), the conclusion follows from \cite[Prop. 4.19]{primer}. For (c), let $\epsilon \in (0,\rho' - \rho - \hatsetd_{\bar\rho}(\epi f, \epi f^\nu))$. First, let $(u,x,\alpha) \in \epi \bar f \cap \ball(\rho)$. Then, $-\eta \leq f(u,x) \leq \rho - a(u) \leq \rho$, which implies that $(u,x,f(u,x)) \in \epi f \cap \ball(\bar\rho)$. We can therefore construct $(\bar u, \bar x, \bar \alpha) \in \epi f^\nu$ such that $\max\{\|\bar u - u\|_2, \|\bar x - x\|_2, |\bar\alpha - f(u,x)| \} \leq \hatsetd_{\bar\rho}(\epi f, \epi f^\nu) + \epsilon$. Let $\delta^\nu = \nsup_{\|u'\|_2\leq \rho'} |a^\nu(u') - a(u')|$. These relations entail that 
\begin{align*}
\bar f^\nu(\bar u, \bar x) & = f^\nu(\bar u, \bar x) + a^\nu(\bar u)\\
& \leq \bar \alpha  + a(u)  + \big|a^\nu(\bar u) - a(\bar u)\big| + \big|a(\bar u) - a(u)\big|\\
& \leq f(u,x) + a(u) + \hatsetd_{\bar\rho}(\epi f, \epi f^\nu) + \epsilon + \delta^\nu + \kappa(\rho')\|\bar u - u\|_2\\
& \leq \alpha + \hatsetd_{\bar\rho}(\epi f, \epi f^\nu) + \epsilon + \delta^\nu + \kappa(\rho')\big(\hatsetd_{\bar\rho}(\epi f, \epi f^\nu) + \epsilon\big),
\end{align*}
where we use the fact that $\|\bar u\|_2 \leq \rho + \hatsetd_{\bar\rho}(\epi f, \epi f^\nu) + \epsilon \leq \rho'$. Thus, $(\bar u, \bar x, \bar f^\nu(\bar u, \bar x))$ is a point in $\epi \bar f^\nu$ no further from $(u, x, \alpha)$ in the chosen norm than 
$\delta^\nu + (1+\kappa(\rho'))(\hatsetd_{\bar\rho}(\epi f, \epi f^\nu) + \epsilon)$. We repeat the argument with the roles of $\bar f$ and $\bar f^\nu$ reversed, and let $\epsilon$ tend to zero to reach the conclusion.\eop

\state Proof of Proposition \ref{prop:compositeAug}. Clearly, $\theta^\nu\|\cdot\|^\alpha\eto \iota_{\{0\}^m}$ and, by Proposition \ref{prop:composite}, $f^\nu\eto f$. Let $x\in \dom \phi$, which means that $G(x) \in \dom h$. Since $X^\nu\sto X$, there is $x^\nu\in X^\nu\to x$. We construct $u^\nu = G(x)-G^\nu(x^\nu)$. Then, $u^\nu\to 0$ and $f^\nu(u^\nu,x^\nu) = g_0^\nu(x^\nu) + h^\nu(G(x)) \to g_0(x) + h(G(x)) = f(0,x)$.
If $\{X, X^\nu, \nu\in\nats\}$ are nonempty, then Theorem \ref{thm:augmentation}(a) yields the conclusion because $\theta^\nu\|u^\nu\|^\alpha\to 0$ by assumption. The case with empty sets is handled trivially.\eop

\state Proof of Proposition \ref{prop:compositequant}. By \cite[Thm. 2.19, Prop. 5.37]{primer}, $u \in \nargmin f_{y}(\,\cdot\,, x)$ if and only if $G(x) + u \in \partial h^{*}(y)$. By assumption, there is $v \in \partial h^{*}(y) \cap K$ so that $\|u\|_2 = \|v -  G(x)\|_2 \leq \hat\rho$.\eop

\noindent {\bf Acknowledgements.} The second author was supported by AFOSR under grant 21RT0484 and ONR under grant N00014-24-1-2492.\\

\bibliographystyle{plain}
\bibliography{refs}

\end{document}